\documentstyle{amsppt}
%%%%%%%%%%%%%%%%%%%%%%%%%%%%%%%%%%%
%----------------------------------
%12point LondMS style
%\magnification = 1200
%\pagewidth{105mm}
%\hcorrection{.7in}
%\vcorrection{.15in}
%%%%%%%%%%%%%%%%%%%%%%%%%%%%%%%%%%
%--------------------------------
%12point Arxiv style
\magnification=1200
\vcorrection{-3truemm}
\pageheight{208truemm}
%-----------------------------------
%-----------------------------------
%Singular variety style
%\magnification=900
%\pagewidth{137truemm}
%\pageheight{195truemm}
%\hoffset .58 in
%\voffset .5 in
%---------------------------------
%---------------------------------
%\pageheight{240mm}
%10pt hcorrection
%\pagewidth{180mm}
%\hcorrection{.55in}
%---------------------------------
%11point hcorrection
%\magnification=1100
%\pagewidth{120mm}
%\hcorrection{.5in}
%\vcorrection{.1in}
%--------------------------------
%-------------------------------
%10point
%\magnification=1000
%\pagewidth{139truemm}
%\pageheight{199truemm}
%\hcorrection{.55in}
%---------------------

%\topmatter
%\title
%{\smc generalised goursat normal form II: \\
%partial prolongations} 
%\endtitle

%\endtopmatter
 
%\vskip 12 pt

\author
{\smc PETER J. VASSILIOU}
\endauthor

\heading \smc A Constructive Generalised Goursat Normal Form\endheading
\vskip 20 pt
\heading\smc Peter J. Vassiliou\endheading
\centerline{School of Mathematics and Statistics}
\centerline{University of Canberra, Canberra, A.C.T., 2601}
\centerline{\tt peterv\@ise.canberra.edu.au}

\vskip 40 pt

\def\a{\alpha}
\def\b{\beta}
\def\d{\delta}
\def\D{\Delta}
\def\e{\varepsilon}

\def\G{\Gamma}

\def\O{\Omega}
\def\f{\phi}

\def\o{\omega}
\def\O{\Omega}
\def\r{\rho}
\def\s{\sigma}
\def\S{\Sigma}

\def\x{\xi}

\def\CE{{\Cal E}}

\def\CV{{\Cal V}}
\def\P{\partial}

\def\P #1{\partial_{#1}}

\def\ch #1{\text{\rm Char}\ #1}

\document

%\baselineskip=2\baselineskip

\heading Abstract\endheading 

We provide necessary and sufficient conditions on the derived type of a vector field distribution $\CV$ in order that it be locally
equivalent to a partial prolongation of the contact distribution $\Cal C^{(1)}_q$, on the $1^{st}$ 
order jet bundle of maps from $\Bbb R$ to
$\Bbb R^q$, $q\geq 1$. This result fully generalises the classical Goursat normal form. Our proof is constructive: it is proven that if $\CV$ is locally equivalent to  a partial prolongation of  
$\Cal C^{(1)}_q$ then the explicit construction of contact coordinates algorithmically depends upon the integration of a sequence of geometrically defined and algorithmically determined integrable Pfaffian systems on the ambient manifold. 

\vskip 30 pt
\item{}{\it 2000 MSC}: {\bf 58J60}
\item{}{\it keywords}: {Goursat normal form, derived type, contact distribution, partial prolongation, local equivalence}

\newpage

\centerline{1. {\it Introduction.}}
\vskip 10 pt
The last two decades have witnessed substantial progress in the creation of geometric characterisations of contact distributions over jet bundles [1,8,4]. The problem of obtaining a simple, invariant characterisation of the corresponding {\it partial prolongations} has apparently remained open, even in the case of partial prolongations of jet bundles of maps from the line. That is, partial prolongations of the contact distribution $\Cal C^{(1)}_q$ over the first order jet bundle of maps $\Bbb R\to \Bbb R^q$, $q\geq 1$. 

Finding a simple characterisation of the partial prolongations turns out to be a problem with many applications. For instance, a version of it has been the focus of considerable attention in nonlinear control theory where the equivalence group is not the full diffeomorphism group of the ambient manifold but rather a subgroup preserving the form of the control system [3]. There are also applications to partial differential equations (PDE) where there is a precise sense in which certain integrable hyperbolic PDE systems can be canonically associated to partial prolongations of $\Cal C^{(1)}_q$ and thereby shown to be explicitly integrable. 

For the purposes of such applications, it is desirable not only to solve the {\it recognition problem} for partial prolongations but additionally, to find a method for explicitly {\it constructing} an equivalence between a differential system and a partial prolongation, whenever one is known to exist.

In this paper we give a solution to both the characterisation and construction problem for partial prolongations of 
$\Cal C^{(1)}_q$. Specifically, we give simple geometric conditions, expressed in terms of the derived type of $\CV\subset TM$ guaranteeing the existence of some local diffeomorphism from $M$ which identifies $\CV$ with some partial prolongation of 
$\Cal C^{(1)}_q$. This result fully generalises the classical Goursat normal form from the theory of exterior differential systems [2, 6]. For the purpose of construction, we identify canonical and algorithmically determined integrable distributions over the ambient manifold whose invariants supply local equivalences. 

\vskip 20 pt

\centerline{2. {\it Preliminaries}}

 \vskip 5 pt
 {\it 2.1. The derived flag.}  Recall that the structure tensor of a sub-bundle $\CV\subset TM$ over manifold $M$ is a 
map
$$
\d_0\:\CV\times\CV\to TM/\CV,
$$
defined by
$$
\d_0(X,Y)={\bold p_0}([X,Y])
$$
where 
$$
{\bold p_0}\:TM\to TM/\CV
$$
is the natural projection.
If $\d_0$ has constant rank on $M$ we define the 
{\it first derived bundle} $\CV^{(1)}$ of $\CV$ as the maximal sub-bundle of $TM$ such that
$$
{\bold p_0}(\CV^{(1)})=\text{Im}~\d_0.
$$
In the constant rank case, we can iterate this procedure, obtaining the sequence of structure tensors
$$
\d_i\:\CV^{(i)}\times \CV^{(i)}\to TM/\CV^{(i)},\ i\geq 1
$$
each assumed to have constant rank on $M$. Define the sub-bundle $\CV^{(i)}\subset TM$,$\ i\geq 1$,
as the maximal sub-bundle of $TM$ satisfying ${\bold p_{i-1}}(\CV^{(i)})=\text{Im}~\d_{i-1}$
where ${\bold p_{i-1}}\:TM\allowbreak\to TM/\CV^{(i-1)}$ is the natural projection.
The sub-bundle $\CV^{(i)}\subset TM$ is the $i^{th}$ {\it derived bundle} of $\CV$. We call a bundle $\CV\subset TM$
{\it regular} if the ranks of all structure tensors are constant on $M$. In this regular case, for dimension reasons there is an integer $k\geq 0$ such that $\CV^{(k+1)}=\CV^{(k)}$. The smallest such integer is called the {\it derived length}
of $\CV$ and we evidently have a flag of sub-bundles
$$
\CV\subset\CV^{(1)}\subset\CV^{(2)}\subset\cdots\subset\CV^{(k)}
$$
called the {\it derived flag} of $\CV$. Plainly, the numbers $\dim\CV^{(i)}$ are diffeomorphism invariants of $\CV$.

\vskip 10 pt
{\it 2.2.  Cauchy bundles.} The Cauchy system or characteristic system $\chi(\CV)$ of $\CV$ is defined by
$$
\chi(\CV)=\big\{X\in\G(\CV)~|~[X,\G(\CV)]\subseteq\G(\CV)\big\}.
$$
Even if $\CV$ is regular, its Cauchy system need not have constant rank on $M$. However, if it does, there is a 
sub-bundle $\text{Char}~\CV\subseteq\CV$, the {\it Cauchy bundle}, whose space of smooth sections is $\chi(\CV)$.
\vskip 5 pt
\noindent{\smc Definition 2.1.} If a regular bundle $\CV\subset TM$ satisfies the condition that $\chi(\CV^{(i)})$ is constant on $M$ for each $i\geq 0$ then $\CV$ is said to be {\it totally regular}. For each $i$, the bundle 
$\text{Char}~\CV^{(i)}$ is defined and called the {\it Cauchy bundle} or {\it characteristic bundle} of $\CV^{(i)}$.
The elements of $\text{Char}~\CV^{(i)}$ are called {\it Cauchy vectors}.
\vskip 5 pt
It is easy to show that any Cauchy bundle is integrable.
\vskip 5 pt
\noindent{\smc Definition 2.2.} Let $\CV\subset TM$ be a totally regular bundle. By the {\it derived type}
of $\CV$ we shall mean the list of sub-bundles \big\{$\CV^{(l)},\text{Char}~\CV^{(l)}\big\}$ for all
$l\geq 0$.
\vskip 10 pt
{\it 2.3. The singular variety of a sub-bundle.} The final tool we require is less well known than the derived type of a bundle, at least in the form in which we shall use it, so we discuss this in more detail in the remainder of this section. 

Let $\CV\subset TM$ be such that its structure tensor $\d_0$ has constant rank on $M$. For the remainder of this section we will not need to retain the subscript 0. For simplicity of notation, we will henceforth cease to distinguish between a bundle and its module of smooth sections. By $\Bbb P\CV$ we denote the projectivisation of $\CV$. 

Fix a basis $X_1,\ldots,X_n$ for $\CV$. If $Z_1,\ldots,Z_s$ complete $X_1,\ldots,X_n$
to a frame on $M$, then there are functions $c_{\a\b}^k$ on $M$, antisymmetric in lower indices, such that
$$
\d(X_\a,X_\b)=c_{\a\b}^k \bar{Z}_k.
$$
where $\bar{Z}_k$ is the coset with representative $Z_k$. The elements of $\Bbb P\CV$ are lines determined by elements 
$e^\a X_\a\in\CV$ and will be denoted by the symbol $E=[e^\a X_\a]$.
If $E$ is a fixed direction in $\CV$ then we may seek directions $F=[f^\a X_\a]$ such that
$$
\d(e^\a X_\a,f^\b X_\b)=0.\eqno(2.1)
$$
Equation (2.1) can be expressed in the form
$$
\sum_{\a<\b=2}^n(e^\a f^\b-e^\b f^\a)c^k_{\a\b}=0,\ 1\leq k\leq s,\eqno(2.2)
$$
viewed as a homogeneous linear system for $F$ and in the alternative form
$$
\s(E)\cdot{\bold f}=0,\eqno(2.3)
$$
where ${\bold f}=(f^1,f^2,\ldots,f^n)$. Letting ${\bold e}=(e^1,e^2,\ldots,e^n)$ we have clearly that 
${\bold f}={\bold e}$ is always a solution of (2.3).
\vskip 5 pt
{\smc Definition 2.3.} The matrix $\s(E)$ in (2.3) determined by equations (2.2) will be called the {\it polar matrix}
of the point $E\in\Bbb P\CV$. Such a line  will be called {\it singular} if its polar matrix has less than generic
rank. The set of all singular lines in $\CV$ will be denoted by the symbol $\text{Sing}(\CV)$. We will denote the set of singular lines in $\CV$ over a point $x\in M$, by $\text{Sing}(\CV)(x)$.
\vskip 5 pt
{\smc Lemma 2.1.} {\it For each $x\in M$, $\text{Sing}(\CV)(x)\subset \Bbb P\CV_x$, is a linear determinantal variety.}
\vskip 2 pt
\noindent {\it Proof.} On the one hand, for each $x\in M$, the polar matrix $\s(E_x)$ has entries which are linear functions
of the affine coordinates in $\Bbb P\CV_x$. On the other, $\text{Sing}(\CV)(x)$ is determined by equating the  minors of 
$\s(E_x)$ to zero.\qed
\vskip 5 pt
{\smc Definition 2.4.} The set $\text{Sing}(\CV)$ of all singular points of $\Bbb P\CV$, will be called the {\it singular variety} of $\CV$.
\vskip 5 pt
We now describe a very useful invariant associated to each point of ${\Bbb P}\CV$. 

The structure tensor $\d$  of $\CV\subset TM$ induces a map
$$
\text{deg}_{\CV}\:\Bbb P\CV\to\Bbb N
$$
well defined by
$$
\text{deg}_{\CV}([X])=\dim~\text{Image}~\D_X,\ [X]\in \Bbb P\CV
$$
where for each $X\in\CV$,
$$
\D_X\:\CV\to TM/\CV.
$$
is defined by 
$$
\D_X(Y)=\d(X,Y),\ \ Y\in\CV.
$$
\vskip 5 pt
{\smc Definition 2.5.} Let $\CV\subset TM$ be a sub-bundle and $[X]\in\Bbb P\CV$. The integer
$\text{deg}_\CV([X])$ will be called the {\it degree} of $[X]$. 
\vskip 5 pt
{\smc Lemma 2.2.} {\it For any $E\in\Bbb P\CV$, $\text{deg}_\CV(E)$ is a diffeomorphism invariant.}
\vskip 5 pt
Note that if $\text{deg}_\CV([X])=0$, then $X$ is a Cauchy vector and hence Lemma 2.2 implies, as a special case, the elementary fact that $\text{Char}~\CV$ is an invariant sub-bundle of $\CV$.

\vskip 5 pt
{\smc Lemma 2.3.} {\it For any point $E\in\Bbb P\CV$, $\text{deg}_\CV(E)=\text{rank}~\s(E)$. }
\vskip 2 pt
\noindent{\it Proof.} In fixed bases for $\CV$ and $TM/\CV$, the polar matrix $\s(E)$ is the matrix of the vector bundle 
morphism $\Delta_X$.\qed

\vskip 5 pt
{\smc Remark 2.1.} Amplifying Lemma 2.2, $\text{Sing}(\CV)$ is a diffeomorphism invariant in the sense that if $\CV_1,\CV_2$ are sub-bundles over
$M_1,M_2$, respectively and there is a diffeomorphism $\phi\: M_1\to M_2$ that identifies them, $\phi_*\CV_1=\CV_2$, then 
$\text{Sing}(\CV_2)$ and $\text{Sing}(\phi_*\CV_1)$ are equivalent as projective subvarieties of $\Bbb P\CV_2$. That is, for each $x\in M_1$, there is an element of the projective linear group $PGL({\CV_2}_{|_{\f(x)}},\Bbb R)$ that identifies 
$\text{Sing}(\CV_2)(\f(x))$ and $\text{Sing}(\phi_*\CV_1)(\f(x))$.
\vskip 5 pt
We hasten to point out that the computation of the singular variety for any given sub-bundle  $\CV\subset TM$ is algorithmic. That is, it involves only differentiation and commutative algebra operations. In practice, one computes the determinantal variety of the generic polar matrix $\s(E)$ as a sub-variety of $\Bbb P\CV$. 
\vskip 10 pt

{\it 2.4. The singular variety in positive degree.} One frequently encounters bundles possessing non-trivial Cauchy bundles and it is then natural and prudent to consider the notion of degree of a line in the 
quotient by $\ch\CV$. Continue with sub-bundle $\CV\subset TM$  and let
$$
\widehat{\pi}\:TM\to TM/\ch\CV:=\widehat{TM}
$$
be the standard projection assigning an element of $TM$ to its coset. The structure tensor $\d$ of $\CV$ drops to the quotient bundle $\widehat{\CV}:=\CV/\ch\CV$, leading to a tensor
$$
\widehat{\d}\:\widehat{\CV}\times\widehat{\CV}\to\widehat{TM}/\widehat{\CV}
$$   
well defined by
$$
\widehat{\d}(\widehat{X},\widehat{Y})=\widehat{\pi}\big([X,Y]\big)\mod\widehat{\CV},\ \ \forall\ 
\widehat{X},\widehat{Y}\in\widehat{\CV}
$$
where $\widehat{X}=\widehat{\pi}(X)$, $\widehat{Y}=\widehat{\pi}(Y)$. By analogy, we call $\widehat{\d}$ the {\it structure tensor} of $\widehat{\CV}$. 
As before, we are able to introduce a map
$$
\text{deg}_{\widehat{\CV}}\:\Bbb P\widehat{\CV}\to\Bbb N
$$
well defined by
$$
\text{deg}_{\widehat{\CV}}([\widehat{X}])=\dim~\text{Image}~\widehat{\D}_{\widehat{X}}
$$
where for $\widehat{X}\in\widehat{\CV}$,
$$
\widehat{\D}_{\widehat{X}}(\widehat{Y})=\widehat{\d}(\widehat{X},\widehat{Y}),\ \widehat{Y}\in\widehat{\CV}.
$$
Here, by $[\widehat{X}]$ we denote the distribution spanned by $X\in\CV$ and $\ch\CV$. 
Again, by analogy we call $\text{deg}_{\widehat{\CV}}([\widehat{X}])$ the {\it degree} of 
$[\widehat{X}]\in\Bbb P\widehat{\CV}$. 
\vskip 5 pt
{\smc Remark 2.2.} All definitions and results of subsection 2.3 hold {\it mutatis mutandis} when the structure tensor $\d$  is replaced by
$\widehat{\d}$. In particular, we have notions of polar matrix and singular variety, as before. Note however, that each point of $\Bbb P\widehat{\CV}$ has degree one or more. 
\vskip 10 pt

{\it 2.5. The resolvent bundle.} We now specialise the discussion of the previous subsections to introduce an associated sub-bundle which is geometrically and algorithmically determined by $\CV$
and that will play a pivitol role in the present work. 

Suppose a totally regular sub-bundle $\CV\subset TM$ of rank $c+q+1$, $q\geq 2, c\geq 0$ is defined on manifold $M$, 
$\dim M=c+2q+1$. Suppose further that $\CV$ satisfies the following additional properties:
\vskip 2 pt
{\itemitem{i)} $\dim\ch\CV=c$, $\CV^{(1)}=TM$
\itemitem{ii)} 
$\widehat{\S}:=\text{Sing}(\widehat{\CV})=\Bbb P\widehat{\Cal B}\approx\Bbb R\Bbb P^{q-1}$ for some rank $q$ sub-bundle $\widehat{\Cal B}\subset\widehat{\CV}$}
   
\vskip 5 pt
{\smc Definition 2.6.} We will call  $(\CV,\Bbb P\widehat{\Cal B})$ (or $(\CV,\widehat{\S})$) satisfying the above conditions a
{\it   Weber structure} on $M$. 
\vskip 5 pt         
Given a  Weber structure
$(\CV,\Bbb P\widehat{\Cal B})$, let $\Cal R_{\widehat{\S}}(\CV)\subset\CV$, denote the largest sub-bundle such that
$$
\widehat{\pi}\big( \Cal R_{\widehat{\S}}(\CV) \big)= \widehat{\Cal B}.\eqno(2.8)
$$
\vskip 5 pt
{\smc Definition 2.7.} The rank $q+c$ bundle $\Cal R_{\widehat{\S} }(\CV)$ will be called the {\it resolvent bundle} associated to the  Weber structure $({\CV},\widehat{\S})$. The bundle $\widehat{\Cal B}$ determined by the singular variety of 
$\widehat{\CV}$ will be called the {\it singular sub-bundle} of the  Weber structure. A  Weber structure
will be said to be {\it integrable} if its resolvent bundle is integrable.
\vskip 5 pt
{\smc Remark 2.3.} Note that an {\it integrable}  Weber structure descends to the quotient of $M$ by the leaves of $\ch\CV$ to be the contact bundle on $J^1(\Bbb R,\Bbb R^q)$. The term `Weber structure' honours Eduard von Weber (1870 - 1934) who, as far as I can tell, was the first to publish a proof of the Goursat normal form [7].
\vskip 5 pt
{\smc Proposition 2.7.} {\it  Let $(\CV,\widehat{\S})$ be a Weber structure on $M$ and 
$\widehat{\Cal B}$ its singular sub-bundle. If $q\geq 3$, then the following are equivalent
\vskip 3 pt
{\itemitem{1)} Its resolvent bundle $ \Cal R_{\widehat{\S}}(\CV) \subset\CV$ is integrable.                   
\itemitem{2)} Each point of $\widehat{\S}=\text{Sing}(\widehat{\CV})$ has degree one
\itemitem{3)} The structure tensor $\widehat{\d}$ of $\widehat{\CV}$ vanishes on $\widehat{\Cal B}$}}
\vskip 2 pt
\noindent{\it Proof.} 
\newline
{\it 1)\ $\Rightarrow$\ 2)}: It follows from (2.8) that
$\ch\CV\subset\Cal R_{\widehat{\S}}(\CV)$. Suppose that $\Cal R_{\widehat{\S}}(\CV)$ is integrable
and spanned by $\xi_1,\ldots,\xi_c,\b_1,\ldots\b_q$. Let $\xi_1,\ldots,\xi_c,\b_1,\ldots\b_q, X$ span
$\CV$. For each $i$ in the range $1\leq i\leq q$ there must be elements $Z_i$ of $TM$ such that
$$
[\b_i,X]\equiv Z_i\mod\CV\eqno(2.9)
$$
and $Z_1\ldots,Z_q$ complete $\xi_1,\ldots,\xi_c,\b_1,\ldots\b_q, X$ to a frame on $M$. If $\widehat{Z}_k$
$:=\widehat{\pi}(Z_k)$ denotes the coset with representative $Z_k$ then by (2.9) we have
$$
\widehat{\d}(\widehat{\b}_k, \widehat{X})=\widehat{Z}_k.\eqno(2.10)
$$  
If $\widehat{E}=[e^0\widehat{X}+e^1\widehat{\b}_1+\cdots +e^q\widehat{\b}_q]$ is an arbitrary element of 
$\Bbb P\widehat{\CV}$ then as in our previous discussion, we can write
$$
\text{deg}_{\widehat{\CV}}(\widehat{E})=\dim~\Big\{~\sum_{k=1}^q(e^0f^k-e^kf^0)\widehat{Z}_k~|~f^k\ \text{arbitrary}\Big\}.
$$
Thus, generically, $\text{deg}_{\widehat{\CV}}(\widehat{E})=q$. To compute the singular lines in $\widehat{\CV}$ we seek 
those points $\widehat{E}\in\Bbb P\widehat{\CV}$ whose polar matrix has rank less than $q$. In this case, since (2.10) is the only non-zero structure, the analogue of equation (2.5) is
$$
e^0f^k-e^kf^0=0,\ 1\leq k \leq q,\eqno(2.11)
$$
leading to polar matrix
$$
\s(\widehat{E})=\left(\matrix -e^1 & e^0 & 0 & 0 & 0 & \cdot & \cdot  & 0\\
                              -e^2 & 0 & e^0 & 0 & 0 & \cdot & \cdot  & 0\\
                              -e^3 & 0 & 0 & e^0 & 0 & \cdot & \cdot   & 0\\
                              \cdot & \cdot &\cdot &\cdot &\cdot & \cdot & \cdot  &\cdot\\
                              \cdot & \cdot &\cdot &\cdot &\cdot & \cdot & \cdot  &\cdot\\
                              -e^q & 0 & 0 & 0 & 0 & \cdot & \cdot  & e^0\\
\endmatrix\right). 
$$
Clearly, $\text{rank}~\s(\widehat{E})<q$ only if $e^0=0$ in which case $\s(\widehat{E})$ has rank 1. Thus, by Lemma 2.3 and Remark 2.2, 
$\text{deg}_{\widehat{\CV}}(\widehat{E})=1$.   
\vskip 5 pt

{\it 2)\ $\Rightarrow$\ 3):} Suppose each point of $\Bbb P\widehat{\Cal B}$ has degree 1 and let $\{\widehat{\b}_k\}$, 
$1\leq k\leq q$ be a basis for $\widehat{\Cal B}$.  Then there is an element $\langle\widehat{Z}\rangle\in
\widehat{TM}/\widehat{\CV}$ such that 
$$
\text{Image}\ \widehat{\D}_{\widehat{\b_1}}=\{\langle\widehat{Z}\rangle\}.\eqno(2.12)
$$ 
 It follows that there are functions $f_{kl}$ 
satisfying $f_{kl}+f_{lk}=0$ such that
$$
\widehat{\d}(\widehat{\b}_k,\widehat{\b}_l)=f_{kl}\langle\widehat{Z}\rangle,\ k,l=1,2\ldots,q.\eqno(2.13)
$$
Suppose $f_{12}\neq 0$. In view of (2.13) and the fact that each $\widehat{\b}_k$ has degree 1, there must be functions 
$a,b$ on $M$ such that
$$
\widehat{\d}(\widehat{\b}_1,\widehat{X})=a\langle\widehat{Z}\rangle,\ 
\widehat{\d}(\widehat{\b}_2,\widehat{X})=b\langle\widehat{Z}\rangle,\eqno(2.14)
$$
where $\widehat{X}$ completes the $\widehat{\b}_k$ to a basis for $\widehat{\CV}$. Equations (2.13) and (2.14) imply that
$$
\aligned
\dim\CV^{(1)}-\dim\CV=&1+\dim~\{[\b_3,X],[\b_4,X],\ldots,[\b_q,X]\}\cr
                     \leq & 1+q-2=q-1<q,
\endaligned
$$
a contradiction. Hence $f_{12}=0$. It follows that $\widehat{\d}(\widehat{\Cal B},\widehat{\Cal B})=0$.
\vskip 5 pt
{\it 3)\ $\Rightarrow$\ 1):} The fact that $\widehat{\d}$ vanishes on $\widehat{\Cal B}$ implies that 
$[\b_k,\b_l]\in\CV$. Since $[\ch\CV,\ch\CV]\allowbreak\subset\CV$ and $[\b_k,\ch\CV]\subset\CV$, resolvent 
bundle $\Cal R_{\widehat{\S}}(\CV)$ satisfies $[\Cal R_{\widehat{\S}}(\CV),\Cal R_{\widehat{\S}}(\CV)]\subset\CV$.
Since $q\geq 3$, a theorem of R. Bryant [1] (see also [4]), 
guarantees that $\Cal R_{\widehat{\S}}(\CV)$ is integrable. \qed
\vskip 5 pt
{\smc Proposition 2.8.} {\it Let $({\CV},\widehat{\S})$ be an integrable Weber structure on $M$. 
Then its resolvent bundle $\Cal R_{\widehat{\S}}(\CV)$ is the unique, maximal,integrable sub-bundle of $\CV$.}
\vskip 2 pt
\noindent {\it Proof.} We first note that for any integrable Weber structure 
$({\CV},\widehat{\S})$, 
if $\widehat{\Cal B}'\subset\widehat{\CV}$ is another rank $q$ sub-bundle such that each 
point in $\Bbb P\widehat{\Cal B}'$  has degree 1 then $\widehat{\Cal B}'=\widehat{\Cal B}$. 

Suppose on the contrary that such a sub-bundle $\widehat{\Cal B}'\subset\widehat{\CV}$,
$\widehat{\Cal B}'\neq \widehat{\Cal B}$ exists. 
It follows that there is an element 
$\widehat{\b}\in\widehat{\Cal B}'$ not in $\widehat{\Cal B}$ such that $\text{deg}_{\widehat{\CV}}([\widehat{\b}])=1$.   
Since the resolvent bundle $\Cal R_{\widehat{\S}}(\CV)$ is integrable, the proof of implication {\it 1)\ $\Rightarrow$\ 2)}
of Proposition 2.7 shows that every point of $\Bbb P\widehat{\Cal B}$ has degree 1.
If $\widehat{\b}_1,\ldots,\widehat{\b}_q$ is a basis for $\widehat{\Cal B}$ then 
$\widehat{\b},\widehat{\b}_1,\ldots,\widehat{\b}_q$ is a basis for $\widehat{\CV}$ and it follows that
$$
\widehat{\d}(\widehat{\b},\widehat{\b}_k)=a_k\widehat{Z},\ 1\leq k\leq q\eqno(2.15)
$$
for some functions $a_k$ on $M$ and element $\widehat{Z}\in\widehat{TM}$.
By the proof of Proposition 2.7 we have $\widehat{\d}(\widehat{\b}_k,\widehat{\b}_l)=0$; this and (2.15) give the contradiction
$\dim\CV^{(1)}-\dim\CV=1$. Thus, $\widehat{\Cal B}'=\widehat{\Cal B}$.

Next, suppose there is a bundle $\CE\subset\CV$ which is integrable, of rank $q+c$ and 
$\CE\neq\Cal R_{\widehat{\S}}(\CV)$. It is easy to show that $\ch\CV\subset\CE$, else $\CV$ is integrable.
Also, one can check that in view of the integrability of $\CE$, each line in the rank $q$ sub-bundle
$\widehat{\CE}:=\widehat{\pi}(\CE)$ has degree 1. It follows that $\widehat{\CE}=\widehat{\Cal B}$. By hypothesis there
is an element $\e\in\CE$ such that $\e\notin\Cal R_{\widehat{\S}}(\CV)$. But we've shown that 
$\widehat{\pi}(\e)\in\widehat{\Cal B}$ and by the surjectivity of $\widehat{\pi}$ there is a 
$\b\in\Cal R_{\widehat{\S}}(\CV)$ such that $\widehat{\pi}(\e)=\widehat{\pi}(\b)$. Hence for some $\xi\in\ch\CV$,
$\e=\b+\xi\in\Cal R_{\widehat{\S}}(\CV)$. This contradiction proves that $\CE=\Cal R_{\widehat{\S}}(\CV)$. \qed

\vskip 5 pt
{\smc Remark 2.4.} Note that checking the integrability of the resolvent bundle is algorithmic.
One computes the singular variety $\text{Sing}(\widehat{\CV})=\Bbb P\widehat{\Cal B}$.
%If $q\geq 3$ one need only check that $\widehat{\d}$ vanishes on $\widehat{\Cal B}$, by Proposition 2.7.
In turn, the singular bundle $\widehat{\Cal B}$ algorithmically determines  $\Cal R_{\widehat{\S}}(\CV)$.
\vskip 5 pt

We conclude this section by mentioning some notation. Firstly, we note that in this paper we work exclusively in the smooth ($C^\infty$) category and all objects and maps will be assumed to be smooth without further notice. Secondly, we will often denote sub-bundles 
$\CV\subset TM$ by a list of  vector fields $X,Y,Z,\ldots$ on $M$ enclosed by braces,
$$
\CV=\big\{X,Y,Z,\ldots\big\}.
$$ 
This will always denote the bundle $\CV$ whose space of sections is the $C^\infty(M)$-module generated by
vector fields $X,Y,Z,\ldots$.

\vskip 20 pt
   
\centerline{3. {\it Partial prolongations and Goursat bundles.}}
\vskip 10 pt
In this section we give a brief coordinate description of partial prolongations and introduce the notion of a Goursat bundle.

The {\it contact distribution} on the $1^{st}$ order jet bundle of maps from 
$\Bbb R\to \Bbb R^q$, $J^{1}(\Bbb R,\Bbb R^q)$, $q\geq 1$ will be denoted by the symbol $\Cal C^{(1)}_q$ and locally expressed in contact coordinates as
$$
\Cal C^{(1)}_q=\big\{\P x+\sum_{j=1}^qz^j_1\P {z^j}, \P {z^j_1}\big\}.\eqno(3.1)
$$
A {\it partial prolongation} of $\Cal C^{(1)}_q$ may be expressed in contact coordinates as a distribution on 
$J^k(\Bbb R,\Bbb R^q)$ of the form
$$
\Cal C(\tau)=\big\{\P x+\sum_{a=1}^t\sum_{j_a=1}^{q_a}\sum_{l_a=1}^{k_a-1} z^{a,j_a}_{l_a+1}\P {z^{a,j_a}_{l_a}}, 
\P {z^{a,j_a}_{k_a}}  \big\}.\eqno(3.2)
$$    
Here and elsewhere in this paper, the symbol $\tau$ denotes the {\it type} of the partial prolongation which is specified by an ordered list of ordered pairs
$$
\tau=\langle q_1,k_1; q_2,k_2;\ldots;q_t,k_t\rangle\eqno(3.3)
$$
indicating that there are $q_a$ variables of order $k_a$ and we use the convention $1\leq k_1<k_2<\cdots <k_t$; the 
$q_i$ are any positive integers. The positive integer $t$ will be called the {\it class} of $\Cal C(\tau)$.

It will sometimes be convenient to express the type of a partial prolongation as an ordered list of $k_t$ non-negative integers 
$$
\tau=\langle \r_1,\r_2,\ldots,\r_{k_t}\rangle\eqno(3.4)
$$
where the $j^{th}$ element $\r_j$ indicates the number of variables of order $j$. In this notation, if each variable $z^j$, 
$1\leq j\leq q$ has the same order $k_j=k$ then the contact distribution
(3.2) has class 1 and its type is $\langle 0,0,\ldots,0,q\rangle$, where $k-1$ zeros precede entry $q$. Such a contact system is a {\it total prolongation} of $\Cal C^{(1)}_q$, denoted, $\Cal C^{(k)}_q$, an instance of a partial prolongation. Note that, in the general case, the derived length of a partial prolongation of type $\tau$ is $k=k_t$.  

For any totally regular sub-bundle $\CV\subset TM$, we have the notion of its derived type. In section 2, we defined the {\it derived type} of a bundle as the list of all derived bundles together with their corresponding Cauchy bundles. We shall frequently abuse notation by using the term `derived type of $\CV$' for  the ordered list of ordered pairs of the form
$$
[[m_0,\chi^0],[m_1,\chi^1],\ldots,[m_k,\chi^k]]
$$
where $m_j$ denotes the rank of the $j^{th}$ derived bundle $\CV^{(j)}$ of $\CV$ and $\chi^j$ denotes the rank of its Cauchy bundle, $\ch\CV^{(j)}$. 

It is important to relate the type of a partial prolongation to its derived type. For this it's convenient to introduce the notions of {\it velocity}, {\it acceleration} and {\it decceleration} of a sub-bundle. 

\vskip 5 pt

{\smc Definition 3.1.} Let $\CV\subset TM$ be a totally regular sub-bundle with derived type 
$$
\d_\CV=[[m_0,\chi^0],[m_1,\chi^1],\ldots,[m_k,\chi^k]].\eqno(3.5)
$$ 
The {\it velocity} of $\CV$ is the ordered list of $k$ integers 
$$
\text{\rm vel}(\CV)=\langle\D_1,\D_2,\ldots,\D_k\rangle, 
$$
where,
$$
\D_j=m_j-m_{j-1},\ 1\leq j\leq k.
$$
The {\it acceleration} of $\CV$ is the ordered list of $k$ integers 
$$
\text{accel}(\CV)=\langle\D^2_2,\D^2_3,\ldots,\D^2_k,\D_k\rangle,
$$
where
$$
\D^2_i=\D_i-\D_{i-1},\ 2\leq i\leq k.
$$
The {\it decceleration} of $\CV$ is the ordered list of $k$ integers 
$$
\text{deccel}(\CV)=\langle -\D^2_2,-\D^2_3,\ldots,-\D^2_k, \D_k\rangle.
$$ 
Note total prolongations $\Cal C^{(k)}_q$ have accelerations of the form
$$
\langle 0,0,\ldots,0,q\rangle,\ q\geq 1,
$$
where there are $k-1$ zeros before the last entry $q$.  The classical Goursat normal form is the case $q=1$ in this family of deccelerations. The main aim of this paper is to generalise this classical result to completely arbitrary decclerations.
\vskip 5 pt
To recognise when a given sub-bundle has
or has not the derived type of a partial prolongation (3.2) we introduce one further canonically associated sub-bundle that plays a crucial role. 

If $\CV$ has derived length $k$ we let
$
\ch\CV^{(j)}_{j-1}
$
denote the intersections $\CV^{(j-1)}\cap\ch\CV^{(j)}$, $1\leq j\leq k-1$. It is easy to see that in every partial prolongation these sub-bundles are non-trivial and integrable.
\vskip 5 pt

{\smc Proposition 3.1.} {\it Let sub-bundle $\CV\subset TM$ be totally regular with velocity and acceleration
$\langle \D_1,\ldots,\D_k\rangle$ and $\langle\D^2_2,\ldots,\D^2_k,\D_k\rangle$, respectively. Then $\CV$ has the derived type
of a partial prolongation $\Cal C(\tau)$ of some type $\tau$ if and only if
$$
\aligned
&m_0=1+P,\  m_1=1+2P,\cr
&m_l=1+(1+l)P+\sum_{j=2}^l (l+1-j)\D^2_j,\ 2\leq l\leq k,\cr
&\chi^j=2m_j-m_{j+1}-1,\ 0\leq j\leq k-1,\cr
&\chi^i_{i-1}=m_{i-1}-1,\ 1\leq i\leq k-1.
\endaligned\eqno(3.6)
$$
where, $P=\sum_{i=1}^k\r_l$, and 
$$ 
m_l=\dim\CV^{(l)},\ \chi^l=\dim\ch\CV^{(l)},\ \chi^l_{l-1}=\dim\ch\CV^{(l)}_{l-1}. 
$$
The type $\tau$ in $\Cal C(\tau)$ is given by the decceleration,
$
\tau=\text{\rm deccel}(\CV).
$}
\vskip 3 pt
\noindent {\it Proof.} The equations $(3.6)_1$ are easily deduced from the local form (3.2). Equation $(3.6)_2$ is obtained by solving the initial value problem
$$
\aligned
m_{j+1}=\D^2_{j+1}+2m_j-m_{j-1},\cr
m_0=1+P,\ m_1=1+2P.
\endaligned
$$
For the final two equations it follows from the local normal form (3.2) that $\chi^j$ satisfies the recurrence relation 
$\chi^{l+1}=\chi^l+\D_{l+1}-\D^2_{l+2}$, $0\leq l\leq k-2$, with $\chi^0=0$; and that $\chi^l_{l-1}$ satisfies 
$\chi^{l+1}_{l}=\chi^l_{l-1}+\D_{l+1}-\D^2_{l+1}$, $1\leq l\leq k-2$ with $\chi^1_0=\D_1$. From these the remaining equations in (3.6) are readily deduced. Finally, observe that the number of variables of order $l$ in (3.2) is given by
$\chi^l-\chi^l_{l-1}=-\D^2_{l+1}$, for $1\leq l\leq k-1$ and that there are $\D_k$ variables of order $k$.  This shows that the type of a partial prolongation is given by $\tau=\text{\rm deccel}(\CV)$. \qed. 
\vskip 5 pt

{\smc Definition 3.2.}  A totally regular sub-bundle $\CV\subset TM$ of derived length $k$ will be called a {\it Goursat bundle of type} $\tau$
if
\item{[i]} $\CV$ has the derived type of a partial prolongation whose type is $\tau=\text{deccel}(\CV)$ 
\item{[ii]} Each intersection $\ch\CV^{(i)}_{i-1}$ is an integrable sub-bundle whose rank, assumed to be constant on $M$,
agrees with the corresponding rank in $\Cal C(\tau)$
\item{[iii]} In case $\D_k>1$, then $\CV^{(k-1)}$ determines an integrable Weber structure on $M$, with singular sub-bundle
$\Bbb P\widehat{\Cal B}\approx\Bbb R\Bbb P^{\D_k-1}$.

\newpage

\centerline{4. {\it Geometric characterisation of partial prolongations.}}
\vskip 10 pt
In this section we establish our main result by giving a constructive solution to the recognition problem for partial prolongations of the contact distribution $\Cal C^{(1)}_q$
on the jet space $J^1(\Bbb R, \Bbb R^q)$.
\vskip 5 pt
{\smc Theorem 4.1.} [Generalised Goursat Normal Form]. {\it Let $\CV\subset TM$ be a Goursat bundle over manifold $M$, of derived length $k>1$ and type 
$\tau=\text{\rm deccel}(\CV)$. 
Then there is an open, dense subset $\hat M\subseteq M$ such that the restriction of $\CV$ to $\hat M$ is locally equivalent to $\Cal C(\tau)$. Conversely any partial prolongation of $\Cal C^{(1)}_q$ is a Goursat bundle.}
\vskip 5 pt
\noindent {\it Proof.} Firstly, by application of the Jacobi identity, it is easily shown that any Goursat bundle over $M$ induces a filtration of $TM$ by integrable sub-bundles. There is a distinction between the cases $\D_k>1$ and $\D_k=1$. We begin with the former case for which we have 
$$
\aligned
\ch\CV^{(1)}_0&\subseteq\ch\CV^{(1)}\subset\cdots\subset\ch\CV^{(j)}_{j-1}\subseteq\ch\CV^{(j)}\subset\cdots\cr
 \cdots&\subset\ch\CV^{(k-1)}_{k-2}\subseteq\ch\CV^{(k-1)}\subset\Cal R_{\widehat{\S}_{k-1}}(\CV^{(k-1)})\subset TM
\endaligned
$$
where $\Cal R_{\widehat{\S}_{k-1}}(\CV^{(k-1)})$ is the resolvent bundle of the integrable Weber structure 
$(\CV^{(k-1)},\allowbreak \widehat{\S}_{k-1})$. Note that for $j$ in the range $1\leq j\leq k-1$, $\ch\CV^{(j)}_{j-1}=\ch\CV^{(j)}$ if and only if $\D^2_{j+1}=0$. 

We will use the convenient notation
$$
\aligned
&\nu_j=\D_j,\ 1\leq j\leq k,\cr
&n_i=\nu_{i+1},\ 0\leq i\leq k-1,\cr  
&N_l=\dim M-m_l,\ 0\leq l\leq k,
\endaligned\eqno(4.1)
$$
in terms of which we have
\vskip 5 pt

{\smc Lemma 4.2.} {\it Let $\CV$ be a Goursat bundle on $M$ of derived length $k$ with type \text{\rm deccel(}$\CV${\rm)}.
Then 
\item{\rm [i]} $N_{j+1}=N_j-n_j,\ 0\leq j\leq k-1$,
\item{\rm [ii]} $\ch\CV^{(l)}$ has $N_l+n_l+1=N_{l-1}+1$ independent invariants, $1\leq l\leq k-1,$
\item{\rm [iii]} $\ch\CV^{(l)}_{l-1}$ has $N_l+\nu_l+1=N_{l-1}+1+\D^2_{l+1}$ independent invariants, $1\leq l\leq k-1.$}

\vskip 3 pt

\noindent{\it Proof.} Elementary calculation using (4.1) and Proposition (3.1).  \qed

\vskip 10 pt

\vskip 5 pt
In section 3 we constructed an explicit bijection between the set of all derived types $\d_{\Cal V}$ of  Goursat bundles and the set of all types $\tau$ for partial prolongations. We now make use of this and Lemma 4.2 in the first of three steps in our proof of Theorem 4.1 by establishing a crucial semi-canonical form.
\vskip 5 pt
{\smc Lemma  4.3.} {\it Let $\CV\subset TM$ be a Goursat bundle over manifold $M$ of derived length $k$ and type
$\tau=\langle \r_1,\r_2,\ldots,\r_k\rangle$=$\text{\rm deccel (}\CV\text{\rm )}$. Then there are local coordinates 
$$
x,x_1,x_2,\ldots,x_{N_0},x_{N_0+1},\ldots,x_{N_0+\nu_1}
$$
and vector fields $X,\P {N_0+1},\ldots,\P {N_0+\nu_1}$, on $M$ which generate $\CV$. Here and elsewhere, $\P {\a}$ denotes the vector field $\P {x_\a}$ and $X$ has the local form
$$
X=\P x +\sum_{\a=1}^{N_0}\x^\a\P {\a}\eqno(4.2)
$$
where the functions $\xi^\a$ satisfy 
$$
\aligned
\frac{\partial\xi^\a}{\partial x_\s}=0,\ N_j+n_j+1\leq \s\leq \dim M-1,\ &1\leq \a\leq N_j,\cr
& 1\leq j\leq k-1
\endaligned\eqno(4.3)
$$
and
$$
\Bigl|\frac{\partial (\xi^{N_{j-1}-n_{j-1}+1},\ldots,\xi^{N_{j-1}})}{\partial (x_{N_{j-1}+1},
\ldots,x_{N_{j-1}+n_{j-1}})}\Bigr|\neq 0,\ 1\leq j\leq k.\eqno(4.4)
$$
Finally, the derived bundles $\CV^{(j)}$ of $\CV$ have local form
$$
\aligned
\CV^{(j)}=\{\P x+\sum_{\a=1}^{N_j}\xi^\a\P \a,\P {N_j+1},&\ldots,\P {N_j+n_j};\cr
                                  \P {N_j+n_j+1},&\ldots,\P {N_j+\nu_j}\}\oplus\ch\CV^{(j)}_{j-1}
\endaligned\eqno(4.5)
$$

for $j$ in the range $1\leq j\leq k-1$}. 
     
\vskip 3 pt
\noindent {\it Proof.} By induction on the order $j$ of the derived bundles $\CV^{(j)}$ of $\CV$. 
Begin by introducing coordinates 
$$
x,x_1,\ldots,x_{N_0},x_{N_0+1},\ldots,x_{{N_0}+\nu_1}
$$
adapted to the filtration
of $TM$ induced by $\CV$. Specifically, for $0\leq j\leq k-2$, functions 
$$
x,x_1,\ldots, x_{N_j}
$$ 
span the invariants of 
$\ch\CV^{(j+1)}_j$ and 
$$
x,x_1,\ldots x_{N_j+\D^2_{j+2}}
$$
span the invariants of $\ch\CV^{(j+1)}$. Finally, 
$$
x,x_1\ldots, x_{N_{k-1}}
$$
span the invariants of the resolvent bundle $\Cal R_{\widehat{\S}_{k-1}}(\CV^{(k-1)})$.

Since $\dim\ch\CV^{(1)}_0=m_0-1$ we deduce that in these adapted coordinates $\CV$ has the local form 
$$
\Big\{\xi^0\P x+\sum_{\a=1}^{N_0}{\xi^\a}\P {\a},\P {N_0+1},\ldots,\P {N_0+\nu_1}\Big\}
$$
for some collection of functions $\x^\a$ on $M$. But $\x^0$ cannot be generically zero, else $\CV$ would have the regular invariant $x$. Consequently, there is a generic subset $\hat M\subseteq M$ in which $\CV$ has the local form
$$
\big\{X=\P x+\sum_{\a=1}^{N_0}\x^\a\P {\a},\P {N_0+1},\ldots,\P {N_0+\nu_1}\big\}.\eqno(4.6)
$$
Hence
$$
\CV^{(1)}=\Big\{X,\sum_{\a=1}^{N_0}\frac{\partial\x^\a}{\partial x_{N_0+j_1}}\P {\a},\P {N_0+j_1}\Big\}_{j_1=1}^{\nu_1}.
$$
Since $\dim\ch\CV^{(2)}_1=m_1-1$ and because $x$ is an invariant of $\ch\CV^{(2)}_1$, we deduce that
$$
\ch\CV^{(2)}_1=\Bigg\{\sum_{\a=1}^{N_0}\frac{\partial\x^\a}{\partial x_{N_0+j_1}}\P {\a},\P {N_0+j_1}\Bigg\}_{j_1=1}^{\nu_1}.
$$
But we recall that the invariants of $\ch\CV^{(2)}_1$ are $x,x_1,\ldots,x_{N_1}$ and hence
$$
\frac{\partial \xi^\a}{\partial x_{N_0+j_1}}=0,\ 1 \leq \a\leq N_1,\ 1\leq j_1\leq \nu_1.\eqno(4.7)
$$
and 
$$
\ch\CV^{(2)}_1=\Bigg\{\sum_{\a=N_1+1}^{N_0}\frac{\partial\x^\a}{\partial x_{N_0+j_1}}\P {\a},\P {N_0+j_1}
\Bigg\}_{j_1=1}^{\nu_1}\eqno(4.8)
$$
in which case we deduce, in view if the rank of $\ch\CV^{(2)}_1$ that
$$
\Bigl|\frac{\partial (\xi^{N_0-n_0+1},\ldots,\xi^{N_{0}})}{\partial (x_{N_{0}+1},
\ldots,x_{N_{0}+n_{0}})}\Bigr|\neq 0.\eqno(4.9)
$$
Consequently,
$$
\aligned
\CV^{(1)}=\{\P x+\sum_{\a=1}^{N_1}\xi^\a\P \a, \P {N_1+1},&\ldots, \P {N_1+n_1};\cr
\P {N_1+n_1+1},&\ldots,\P {N_1+\nu_1}\}\oplus\ch\CV^{(1)}_0
\endaligned\eqno(4.10)
$$
since, by construction, the vector fields $\P {N_0+j_1}$, $1\leq j_1\leq \nu_1$ form a basis for $\ch\CV^{(1)}_0$. By virtue of our coordinate system, we have
$$
\ch\CV^{(1)}=\{\P {N_1+n_1+1},\ldots,\P {N_1+\nu_1}\}\oplus\ch\CV^{(1)}_0
$$
and consequently
$$
\frac{\partial\xi^\a}{\partial x_{N_1+i_1}}=0,\ 1\leq \a\leq N_1,\ n_1+1\leq i_1\leq \nu_1.\eqno(4.11)
$$
Combining (4.7) and (4.11) we have
$$
\frac{\partial \xi^\a}{\partial x_\s}=0,\ N_1+n_1+1\leq\s\leq\dim M-1,\ 1\leq \a\leq N_1.\eqno(4.12)
$$
Equations (4.12) and (4.9) establish Lemma 4.3 in the case $j=1$.
\vskip 5 pt
Suppose now that (4.3)-(4.5) hold for all $j=1,2,\ldots,l-1<k-2$. Consequently, we have, in particular,
$$
\aligned
\frac{\partial \xi^\a}{\partial x_\s}=0,\ N_{l-1}+n_{l-1}+1&\leq \s\leq \dim M-1,\cr
&\ \ \ \ \ \ \ \ \ \ \ \ \ \ \ \ \ \ \ \ 1\leq\a\leq N_{l-1}
\endaligned\eqno(4.13)
$$
and
$$
\aligned
\CV^{(l-1)}=\Big\{\P x+\sum_{\a=1}^{N_{l-1}}\xi^\a\P \a, \P {N_{l-1}+1},&\ldots, \P {N_{l-1}+n_{l-1}};\cr
\P {N_{l-1}+n_{l-1}+1},&\ldots,\P {N_{l-1}+\nu_{l-1}}\Big\}\oplus\ch\CV^{(l-1)}_{l-2}.
\endaligned\eqno(4.14)
$$
Because of (4.14), we have
$$
\aligned
\CV^{(l)}=\Big\{\P x +\sum_{\a=1}^{N_{l-1}}\xi^\a\P \a, 
\sum_{\a=1}^{N_{l-1}}\frac{\partial\xi^\a}{\partial x_{N_{l-1}+i_{l-1}}}\P \a,&\P {N_{l-1}+i_{l-1}},
\Big\}_{i_{l-1}=1}^{n_{l-1}}\cr
&\oplus\ch\CV^{(l-1)}.
\endaligned\eqno(4.15)
$$
By Proposition 3.1, $\ch\CV^{(l+1)}_l$ has codimension 1 in $\CV^{(l)}$ and invariant $x$ and so
$$
\aligned
\ch\CV^{(l+1)}_l=&\Big\{\sum_{\a=1}^{N_{l-1}}\frac{\partial\xi^\a}{\partial x_{N_{l-1}+i_{l-1}}}\P \a,
\P {N_{l-1}+1},\ldots,\cr
&\ \ \ \ \ \ \ \ \ \ \ \ \ \ \ \ \ \ \ \ \ \ \ \ \ \ \ \ldots,\P {N_{l-1}+n_{l-1}}\Big\}\oplus\ch\CV^{(l-1)}.
\endaligned\eqno(4.16)
$$
By Lemma 4.2, $\ch\CV^{(l+1)}_l$ has $N_{l+1}+\nu_{l+1}+1=N_{l-1}-n_{l-1}+1$ invariants
$$
x,x_1,\ldots,x_{N_{l-1}-n_{l-1}}.
$$
Consequently, 
$$
\frac{\partial \xi^\a}{\partial x_{N_{l-1}+i_{l-1}}}=0,\ 1\leq\a\leq N_{l-1}-n_{l-1},\ 1\leq i_{l-1}\leq n_{l-1}\eqno(4.17)
$$
and so the sub-bundle
$$
\Big\{\sum_{\a=1}^{N_{l-1}}\frac{\partial\xi^\a}{\partial x_{N_{l-1}+i_{l-1}}}\P \a\Big\}_{i_{l-1}=1}^{n_{l-1}}
\subset\ch\CV^{(l+1)}_l
$$
has the form
$$
\Big\{\sum_{\a=N_{l-1}-n_{l-1}+1}^{N_{l-1}}
\frac{\partial\xi^\a}{\partial x_{N_{l-1}+i_{l-1}}}\P \a\Big\}_{i_{l-1}=1}^{n_{l-1}}.\eqno(4.18)
$$
But as $n_{l-1}=\nu_l=m_l-m_{l-1}$, it follows that (4.18) has full rank $n_{l-1}$ and hence
$$
\Bigl|\frac{\partial (\xi^{N_{l-1}-n_{l-1}+1},\ldots,\xi^{N_{l-1}})}{\partial (x_{N_{l-1}+1},
\ldots,x_{N_{l-1}+n_{l-1}})}\Bigr|\neq 0.\eqno(4.19)
$$
Because of (4.19), we can express $\CV^{(l)}$ in (4.15) in the form
$$
\aligned
\CV^{(l)}=&\Big\{\P x+\sum_{\a=1}^{N_{l-1}-n_{l-1}}\xi^\a\P \a,\P {N_{l-1}-n_{l-1}+1},\ldots,\P {N_{l-1}};\cr
&\ \ \ \ \ \ \ \ \ \ \ \ \ \ \ \ \ \ \ \ \P {N_{l-1}+1},\ldots,\P {N_{l-1}+n_{l-1}}\Big\}\oplus\ch\CV^{(l-1)}.
\endaligned
$$
However, because of (4.17) and the fact that $n_{l-1}=\nu_l$, we can express $\CV^{(l)}$ in the form
$$
\aligned
\CV^{(l)}=&\Big\{\P x + \sum_{\a=1}^{N_l}\xi^\a\P \a,\P {N_l+1},\ldots,\P {N_l+\nu_l}\Big\}\oplus\ch\CV^{(l)}_{l-1}\cr
         =&\Big\{\P x + \sum_{\a=1}^{N_l}\xi^\a\P \a,\P {N_l+1},\ldots,\P {N_l+n_l};\cr
          &\ \ \ \ \ \ \ \ \ \ \ \ \ \ \ \ \ \ \ \ \ \ \ \ \P {N_l+n_l+1},\ldots,\P {\nu_l}\Big\}\oplus\ch\CV^{(l)}_{l-1}. 
\endaligned\eqno(4.20)
$$
Since $\ch\CV^{(l)}=\{\P {N_l+n_l+1},\ldots,\P {N_l+\nu_l}\}\oplus\ch\CV^{(l)}_{l-1}$, we have that
$$
\frac{\partial\xi^\a}{\partial x_{N_l+j_l}}=0,\ 1\leq\a\leq N_l,\ n_l+1\leq j_l\leq \nu_l.\eqno(4.21)
$$
Combining (4.21), (4.17) with the inductive hypothesis gives that
$$
\frac{\partial\xi^\a}{\partial x_\s}=0,\ N_l+n_l+1\leq \s\leq \dim M-1,\ 1\leq \a\leq N_l.
$$
This equation together with (4.19) and (4.20) constitute the predicate for $j=l$ and Lemma 4.3 is proved up to $j=k-2$.
Thus, we can assert that (4.3)-(4.5) hold for $j=k-2$ and hence 
$$
\aligned
\CV^{(k-1)}=\Big\{\P x +&\sum_{\a=1}^{N_{k-2}}\xi^\a\P \a, 
\sum_{\a=1}^{N_{k-2}}\frac{\partial\xi^\a}{\partial x_{N_{k-2}+i_{k-2}}}\P \a,\P {N_{k-2}+1},\ldots,\cr
&\ldots,\P {N_{k-2}+n_{k-2}}\Big\}\oplus\ch\CV^{(k-2)}
\endaligned
$$
where $1\leq i_{k-2}\leq n_{k-2}$. Property [iii] of the definition of a Goursat bundle implies that its resolvent bundle 
$\Cal R_{\widehat{\S}_{k-1}}(\CV^{(k-1)})$  has co-dimension 1 in $\CV^{(k-1)}$. Since $x$ is one of its invariants, we must have that
$$
\Big\{\sum_{\a=1}^{N_{k-2}}\frac{\partial\xi^\a}{\partial x_{N_{k-2}+i_{k-2}}}\P \a\Big\}_{i_{k-2}=1}^{n_{k-2}}\subset
\Cal R_{\widehat{\S}_{k-1}}(\CV^{(k-1)}).\eqno(4.22)
$$
But the resolvent bundle has invariants $x,x_1,\ldots,x_{N_{k-1}}$ and since $N_{k-1}=N_{k-2}-n_{k-2}$, we deduce from (4.22) that
$$
\frac{\partial \x^\a}{\partial x_{N_{k-2}+i_{k-2}}}=0,\ 1\leq \a\leq N_{k-1},\ 1\leq i_{k-2}\leq n_{k-2},\eqno(4.23)
$$
in which case the lower limit in the left hand side of (4.22) becomes $\a=N_{k-2}+n_{k-2}+1$. But 
the bundle on the left hand side of (4.22) must have rank $m_{k-1}-m_{k-2}=\nu_{k-1}=n_{k-2}$ and hence we deduce
that (4.4) holds for $j=k-1$. This fact and Proposition 3.1 enables us to express $\CV^{(k-1)}$ in the form
$$
\aligned
\CV^{(k-1)}=&\Big\{\P x +\sum_{\a=1}^{N_{k-2}-n_{k-2}}\xi^\a\P \a,
\P {N_{k-2}-n_{k-2}+1},\ldots,\P {N_{k-2}};\cr
&\ \ \ \ \ \ \ \ \ \ \P {N_{k-2}+1},\ldots,\P {N_{k-2}+n_{k-2}}\Big\}\oplus\ch\CV^{(k-2)}\cr
         =&\Big\{\P x +\sum_{\a=1}^{N_{k-1}}\xi^\a\P \a,\P {N_{k-1}+1},\ldots,\P {N_{k-1}+n_{k-2}}\Big\}\cr
          &\ \ \ \ \ \ \ \ \ \ \ \ \ \ \ \ \ \ \ \ \ \ \ \ \ \ \ \ \ \ \ \ \ \ \ \ \ \ \ \ \ \ \ \ \ \ \ 
\oplus\ch\CV^{(k-1)}_{k-2}. 
\endaligned\eqno(4.24)
$$
Because $n_{k-2}=\nu_{k-1}$, by virtue of our coordinate system, we have that
$$
\ch\CV^{(k-1)}=\{\P {N_{k-1}+n_{k-1}+1},\ldots,\P {N_{k-1}+\nu_{k-1}}\}\oplus\ch\CV^{(k-1)}_{k-2}
$$ 
and hence from (4.24) we deduce that
$$
\frac{\partial\xi^\a}{\partial x_{N_{k-1}+j_{k-1}}}=0,
1\leq\a\leq N_{k-1},\ n_{k-1}+1\leq j_{k-1}\leq \nu_{k-1}.\eqno(4.25)
$$
Combining (4.23), (4.25) with (4.3) for $j=k-2$ permits us to conclude that (4.3) holds for $j=k-1$. From this fact and
(4.24) we deduce that
$$
\aligned
\CV^{(k)}=\Big\{\P x +&\sum_{\a=1}^{N_{k-1}}\xi^\a\P \a, 
\sum_{\a=1}^{N_{k-1}}\frac{\partial\xi^\a}{\partial x_{N_{k-1}+i_{k-1}}}\P \a,\P {N_{k-1}+1},\ldots,\cr
&\ldots,\P {N_{k-1}+n_{k-1}}\Big\}\oplus\ch\CV^{(k-1)}
\endaligned\eqno(4.26)
$$
with $1\leq i_{k-1}\leq n_{k-1}$. But as $n_{k-1}=\nu_k=m_k-m_{k-1}$, the bundle
$$
\Big\{\sum_{\a=1}^{N_{k-1}}\frac{\partial\xi^\a}{\partial x_{N_{k-1}+i_{k-1}}}\P \a\Big\}_{i_{k-1}=1}^{n_{k-1}}\eqno(4.27)
$$
must have rank $n_{k-1}$. But note that, in fact $n_{k-1}=N_{k-1}$ and this implies that (4.4) holds for $j=k$.
Lemma 4.3 is now proved in case $\D_k>1$. 

Now suppose $\D_k=1$. In this case the resolvent bundle is not defined and the second last bundle in the filtration induced by $\CV$ is $\ch\CV^{(k-1)}$, in which case Lemma 4.3 holds for all orders up to and including $j=k-2$. Hence, from (4.5) with $j=k-2$ we deduce that
$$
\CV^{(k-1)}=\big\{\P x+\sum_{\a=1}^{N_{k-2}}\x^\a\P \a,X_{i_{k-2}},     
\P {N_{k-2}+i_{k-2}}\big\}_{i_{k-2}=1}^{n_{k-2}}\oplus\ch\CV^{(k-2)}
$$
where
$$
X_{i_{k-2}}=\sum_{\a=1}^{N_{k-2}}\frac{\partial \x^\a}{\partial x_{N_{k-2}+i_{k-2}}}\P \a.
$$
It is easy to deduce from Proposition 3.1 that if $\D_k=1$ then $\chi^{(k-1)}=m_{k-1}-2$. This together with the fact that 
$x$ is an invariant of $\ch\CV^{(k-1)}$ implies that the later is a codimension 1 subspace of 
$$
\Pi^k:=\{X_{i_{k-2}},\P {N_{k-2}+i_{k-2}}\}_{i_{k-2}=1}^{n_{k-2}}\oplus\ch\CV^{(k-2)}
$$
Furthermore, we have 
$$
\ch\CV^{(k-1)}_{k-2}=\{\P {N_{k-2}+i_{k-2}}\}_{i_{k-2}=1}^{n_{k-2}}\oplus\ch\CV^{(k-2)}.
$$
Suppose
$$
Y_l=\eta^{i_{k-2}}_lX_{i_{k-2}}\mod\ch\CV^{(k-1)}_{k-2},\ 1\leq l\leq n_{k-2}-1,
$$
is a basis for $\ch\CV^{(k-1)}/\ch\CV^{(k-1)}_{k-2}$ and $Y\in\Pi^k$ does not belong to $\ch\CV^{(k-1)}$. It follows that
$$
[Y,Y_l]=aX+b^r_lY_r\mod\ch\CV^{(k-1)}_{k-2}
$$
for some functions $a,b^r_l$. As $x$ is an invariant of $[Y,Y_l]$ we must have $a=0$, showing that $\Pi^k$ is integrable and has rank $m_{k-1}-1$. 

By similar reasoning, we easily deduce that 
$$
\Pi^{k+1}:=\text{\rm Im}~\d_{k-1}\oplus\ch\CV^{(k-1)}
$$
is integrable and has dimension $m_k-1$.

We have shown that if $\D_k=1$, then in coordinates adapted to the filtration up to $\ch\CV^{(k-1)}$, there are locally defined integrable bundles $\Pi^k, \Pi^{k+1}$ of ranks $m_{k-1}-1, m_k-1$, respectively which refine the filtration to be
$$
\aligned
\ch\CV^{(1)}_0&\subseteq\ch\CV^{(1)}\subset\cdots\subset\ch\CV^{(j)}_{j-1}\subseteq\ch\CV^{(j)}\subset\cdots\cr
 \cdots&\subset\ch\CV^{(k-1)}_{k-2}\subseteq\ch\CV^{(k-1)}\subset\Pi^k\subset\Pi^{k+1}\subset TM.
\endaligned\eqno(4.28)
$$
We may now choose a coordinate system adapted to filtration (4.28). In particular, $\Pi^k$ has invariants
$$
x,x_1,\ldots,x_{N_{k-1}}.
$$
But $1=m_k-m_{k-1}=\dim M-m_{k-1}-(\dim M-m_k)=N_{k-1}$, so we have that $\Pi^k$ has 2 invariants, $x,x_1$. Furthermore 
$\ch\CV^{(k-1)}$ has invariants $x,x_1,\ldots,\allowbreak x_{N_{k-2}+\D^2_k}$. But $N_{k-1}+\D^2_k=N_{k-1}+1=2$. We deduce that
$$
\Pi^k:=\Big\{\sum_{\a=N_{k-1}+1}^{N_{k-2}}\frac{\partial \x^\a}{\partial x_{N_{k-2}+i_{k-2}}}\P \a
,\P {N_{k-2}+i_{k-2}}\Big\}_{i_{k-2}=1}^{n_{k-2}}\oplus\ch\CV^{(k-2)}.
$$
Now, $N_{k-2}-N_{k-1}=N_{k-2}-(N_{k-2}-n_{k-2})=n_{k-2}$. Since $\Pi^k$ has rank $m_{k-1}-1$ we deduce that (4.4) holds in the case $j=k-1$ and consequently
$$
\aligned
\CV^{(k-1)}=\{\P x +\x^1\P {x_1},\P {N_{k-1}+1}&,\ldots,\P {N_{k-2}},\P {N_{k-2}+1},\cr
&\dots,\P {N_{k-2}+n_{k-2}}\}\oplus\ch\CV^{(k-2)}.
\endaligned
$$
Since $\ch\CV^{(k-1)}$ has invariants $\{x,x_1,\ldots,x_{N_{k-1}+1}\}=\{x,x_1,x_2\}$, it follows that 
$$
\ch\CV^{(k-1)}=\{\P {x_3},\P {x_4},\ldots\}
$$
and hence
$$
\frac{\partial \xi^1}{\partial x_\s}=0,\ 3\leq \s\leq \dim M-1.\eqno(4.29)
$$
Since $\CV^{(k-1)}$ is not integrable, it follows that 
$$
\frac{\partial \x^1}{\partial x_2}\neq 0.\eqno(4.30)
$$
Equation (4.29) is the predicate (4.3) in the case $j=k-1$ and (4.30) is the predicate (4.4) in the case $j=k$. Lemma 4.3 is now proved. \qed  
\vskip 5 pt

With this control over the local form of $\CV$, we are ready to construct contact coordinates. We define vector-valued functions $F_l$, $0\leq l\leq k$, inductively, by
$$
F_0=\left(\matrix x_1 \\ x_2\\ \cdot \\ \cdot \\ x_{N_{k-1}}\endmatrix\right),\ 
F_{p+1}=\left(\matrix XF_p \\ x_{N_{k-(p+1)}+n_{k-(p+1)}+1}\\ \cdot \\ \cdot \\ 
x_{N_{k-(p-1)}+\nu_{k-(p-1)}}\endmatrix\right),\ 0\leq p\leq k-1.\eqno(4.31)
$$
\vskip 5 pt
{\smc Definition 4.1.} We will say that a locally defined real-valued function $\phi$ on $M$ has {\it rank}\  $\sigma$ if $\phi$ depends at most upon the coordinates $x,x_1,\ldots,x_\sigma$ for some $\sigma\leq \dim M$. 
\vskip 5 pt
{\smc Lemma 4.4.} {\it Let $\CV\subset TM $ be a Goursat bundle over manifold $M$ of derived length $k$ and type 
$\text{\rm deccel}(\CV)$. Suppose 
$$
\CV=\{X,\P {N_0+1},\ldots,\P {N_0+\nu_1}\}
$$
is its semi-canonical form. Then the functions $F_l$, defined by (4.31) have the form
$$
F_l=M_l\left(\matrix \x^{N_{k-(l-1)}+1} \\ \cdot \\ \cdot \\ \x^{N_{k-l}} \\ x_{N_{k-l}+n_{k-l}+1} \\
        \cdot \\ \cdot \\   x_{N_{k-l}+\nu_{k-l}}\endmatrix\right)
+f_{[N_{k-(l-1)}+n_{k-(l-1)}]}\eqno(4.32)
$$
for $2\leq l\leq k-1$, where,
$$
M_l=\Bigg(\matrix M_{l-1}\widehat{\G}_{l-1} & 0_{\nu_{k-(l-1)}\times\r_{k-l}}\\
          0_{\r_{k-l}\times\nu_{k-(l-1)}} & I_{\r_{k-l}}\endmatrix\Bigg),\ 2\leq l\leq k-1,\eqno(4.33)
$$
$$ 
\widehat{\G}_{l-1}=\Bigg(\matrix \G_{l-1} & 0_{n_{k-(l-1)}\times\r_{k-l}}\\
          0_{\r_{k-l}\times n_{k-(l-1)}} & I_{\r_{k-l}}\endmatrix\Bigg),\ 2\leq l\leq k.\eqno(4.34)        
$$
The symbol $f_{[\s]}$ indicates functions of rank $\s$ and 
$\r_a=\nu_a-n_a$ is the $a^{th}$ entry in the decceleration vector of $\CV$. The entries
of matrix $M_l$ are functions of rank $N_{k-(l-1)}+n_{k-(l-1)}$
and $\G_{k-(j-1)}$ is the $j^{th}$ matrix in the sequence (4.4).}

\vskip 3 pt

\noindent{\it Proof.} By induction over $l$. The symbol $0_{a\times b}$ in (4.33) and (4.34) denotes a zero matrix and $I_c$ denotes the 
$c\times c$ identity matrix. By (4.31)
$$
F_1=\left(\matrix\x^1\\ \cdot\\ \cdot\\ \x^{N_{k-1}} \\ x_{N_{k-1}+n_{k-1}+1} \\ \cdot\\ \cdot \\ x_{N_{k-1}+\nu_{k-1}}\endmatrix\right),\ 
F_2=\left(\matrix X\x^1\\ \cdot\\ \cdot\\ X\x^{N_{k-1}} \\ \x_{N_{k-1}+n_{k-1}+1} \\ \cdot\\ \cdot\\
\x_{N_{k-1}+\nu_{k-1}} \\ x_{N_{k-2}+n_{k-2}+1} \\ \cdot \\ \cdot \\ x_{N_{k-2}+\nu_{k-2}} \endmatrix\right).\eqno(4.35)
$$
By Lemma 4.3, equation (4.3), the functions $\xi^1,\ldots,x^{N_{k-1}}$ have rank $N_{k-1}+n_{k-1}$. So if 
$h\in\{1,\ldots,N_{k-1}\}$ then
$$
\aligned
X\x^h&=\frac{\partial \x^h}{\partial x}+\sum_{\a=1}^{N_{k-1}+n_{k-1}}\x^\a\frac{\partial\x^h}{\partial x_\a}\cr
&=\frac{\partial \x^h}{\partial x}+\sum_{\a=1}^{N_{k-1}+n_{k-1}}\x^\a\frac{\partial\x^h}{\partial x_\a}\cr
&=\frac{\partial \x^h}{\partial x}+\sum_{\a=1}^{N_{k-1}}\x^\a\frac{\partial\x^h}{\partial x_\a}+
\sum_{\a=N_{k-1}+1}^{N_{k-1}+n_{k-1}}\x^\a\frac{\partial\x^h}{\partial x_\a}. 
\endaligned
$$
Consequently,
$$
X\vec{\x}_1=f_{[N_{k-1}+n_{k-1}]}+\G_1\vec{\x}_{N_{k-1}}
$$
where
$\G_1$ is the matrix in (4.4) with $j=k$ and
$$
\aligned
&\vec{\x}_1=\left(\matrix \x^1 & \cdot & \cdot & \x^{N_{k-1}}\endmatrix\right)^{T},\cr
&\vec{\x}_{N_{k-1}}=\left(\matrix \x^{N_{k-1}+1} & \cdot & \cdot & \x^{N_{k-1}+n_{k-1}}\endmatrix\right)^{T}.
\endaligned
$$
Thus,
$$
F_2=\left(\matrix \G_1 \vec{\x}_{N_{k-1}} \\ \x^{N_{k-1}+n_{k-1}+1} \\ \cdot \\ 
\cdot \\ \x^{N_{k-1}+\nu_{k-1}} \\ x_{N_{k-2}+n_{k-2}+1}\\ \cdot \\ \cdot \\ x_{N_{k-2}+\nu_{k-2}}\endmatrix\right)
+f_{[N_{k-1}+n_{k-1}]}
$$
and hence
$$
F_2=\left(\matrix \widehat{\G}_1 {\bold x}_{N_{k-1}}  \\ x_{N_{k-2}+n_{k-2}+1}\\ \cdot \\ \cdot \\ 
x_{N_{k-2}+\nu_{k-2}}\endmatrix\right)
+f_{[N_{k-1}+n_{k-1}]}
$$
where 
$$
\widehat{\G}_1=\left(\matrix \G_1 & 0_{n_{k-1}\times\r_{k-1}}\\
                              0_{n_{k-1}\times\r_{k-1}} & I_{\r_{k-1}}\endmatrix\right),
{\bold x}_{N_{k-1}}=\left(\matrix \x^{N_{k-1}+1} \\ \cdot \\ \cdot \\ \x^{N_{k-2}}\endmatrix\right).
$$
Consequently,
$$
F_2=M_2\left(\matrix {\bold x}_{N_{k-1}} \\ x_{N_{k-2}+n_{k-2}+1}\\ \cdot \\ \cdot \\ x_{N_{k-2}+\nu_{k-2}}\endmatrix\right)
+f_{[N_{k-1}+n_{k-1}]}\eqno(4.36)
$$
where
$$
M_2=\left(\matrix M_1\widehat{\G}_1 & 0_{\nu_{k-1}\times\r_{k-2}}\\
                              0_{\nu_{k-1}\times\r_{k-2}} & I_{\r_{k-2}}\endmatrix\right).
$$
and $M_1=I_{\nu_{k-1}}$. Equation (4.36) is the predicate in case $l=2$ since the entries of matrix $M_2$ are functions of rank $N_{k-1}+n_{k-1}$. 

\vskip 5 pt

Suppose the inductive hypothesis is satisfied for all $s\leq l$ for some $l<k-1$. Then
$$
F_l=M_l\left(\matrix {\bold x}_{N_{k-(l-1)}} \\ x_{N_{k-(l-1)}+n_{k-(l-1)}+1} \\ \cdot \\ \cdot \\ x_{N_{k-l}+\nu_{k-l}}
\endmatrix\right)+f_{[N_{k-(l-1)}+n_{k-(l-1)}]}\eqno(4.37)
$$
where  
$$
{\bold x}_{N_{k-(l-1)}}=\left(\matrix \x^{N_{k-(l-1)}+1} & \cdot & \cdot & \x^{N_{k-l}}\endmatrix\right)^{T}
$$
and the entries of matrix $M_l$ are functions of rank $N_{k-(l-1)}+n_{k-(l-1)}$.
The following is easily established, using equation (4.3). 
\vskip 5 pt

{\it Claim :} 
{\it If a function $\phi$ on $M$ has rank $N_{k-(l-1)}+n_{k-(l-1)}$
then $X\phi$ has rank $N_{k-l}+n_{k-l}$}. 
\vskip 5 pt
Using this Claim and equation (4.37) we deduce that
$$
F_{l+1}=\left(\matrix M_lX\Xi_l \\ x_{N_{k-(l+1)}+n_{k-(l+1)}+1} \\ \cdot \\ \cdot \\ x_{N_{k-(l+1)}+\nu_{k-(l+1)}}\endmatrix\right)+f_{[N_{k-(l-1)}+n_{k-(l-1)}]}\eqno(4.38)
$$
where  
$$
\Xi_l=\left(\matrix  {\bold x}_{N_{k-(l-1)}} \\ x_{N_{k-l}+n_{k-l}+1} \\ \cdot \\ \cdot \\ x_{N_{k-l}+\nu_{k-l}}
\endmatrix\right).
$$
By Lemma 4.3, equation (4.3) and the Claim, we deduce that
$$
X{\bold x}_{N_{k-(l-1)}}=f_{[N_{k-l}+n_{k-l}]}+\G_l\vec{\x}_{N_{k-l}}
$$
where
$$
\vec{\x}_{N_{k-l}}=\left(\matrix \x^{N_{k-l}+1} & \cdot & \cdot & \x^{N_{k-l}+n_{k-l}}\endmatrix\right)^{T}.
$$
We compute that
$$
X\Xi_l
=\widehat{\G}_l{\bold x}_{N_{k-l}}+f_{[N_{k-l}+n_{k-l}]}\eqno(4.39)
$$
where
$$
\widehat{\G}_l=\left(\matrix \G_l & 0_{\nu_{k-l}\times\r_{k-l}}\\
                     0_{\r_{k-l}\times\nu_{k-l}} & I_{\r_{k-l} }\endmatrix\right).
$$
By (4.38) and (4.39) we have
$$
\aligned
F_{l+1}=&\left(\matrix M_l\widehat{\G}_l{\bold x}_{N_{k-l}} \\ x_{N_{k-(l+1)}+n_{k-(l+1)}+1} \\ \cdot \\ \cdot \\ x_{N_{k-(l+1)}+\nu_{k-(l+1)}}\endmatrix\right)+f_{[N_{k-l}+n_{k-l}]}\cr
       =&M_{l+1}\left(\matrix {\bold x}_{N_{k-l}} \\ x_{N_{k-(l+1)}+n_{k-(l+1)}+1} \\ \cdot \\ \cdot \\ x_{N_{k-(l+1)}+\nu_{k-(l+1)}}\endmatrix\right)+f_{[N_{k-l}+n_{k-l}]}
\endaligned
$$
where
$$
M_{l+1}=\left(\matrix M_l\widehat{\G}_l & 0_{\nu_{k-l}\times \r_{k-(l+1)}}\\
                          0_{\r_{k-(l+1)}\times \nu_{k-l}} & I_{\r_{k-(l+1)}}\endmatrix\right).
$$
This proves Lemma 4.4. \qed

\vskip 5 pt

We conclude the proof of Theorem 4.1 by showing that the components of the vector-valued functions $F_l$, $0\leq l\leq k-1$
constructed in Lemma 4.4, together with the functions $x_{N_0+1},\ldots,x_{N_0+\nu_1}$ determine a local coordinate system on $M$. 
\vskip 5 pt
{\smc Lemma 4.5.} {\it Let $\CV\subset TM$ be a Goursat bundle over $M$ of derived length $k$, type $\text{\rm deccel}(\CV)$ and $F_l$, $0\leq l\leq k$ the collection of vector valued functions constructed in Lemma 4.4. Let $x_{N_0+1},\ldots,x_{N_0+\nu_1}$ be the functions constructed in Lemma 4.3. Define differential forms $\omega_j$ for $0\leq j\leq k-1$ inductively by
$$
\omega_0=dF^1_0\wedge dF^2_0\wedge\cdots\wedge dF^{\nu_k}_0,\ \omega_j=\omega_{j-1}\wedge dF^1_j\wedge dF^2_j
\wedge\cdots\wedge dF^{\nu_{k-j}}_j.\eqno(4.40)
$$
Then 
$$
\omega_j\approx(\det\G_1)^j(\det\G_2)^{j-1}\cdots (\det\G_j)\hskip 0.75 pt\text{\rm vol}_{N_{k-(j+1)}},\ 
0\leq j\leq k-1\eqno(4.41)
$$
where $\text{\rm vol}_\sigma=dx\wedge dx_1\wedge \cdots\wedge dx_\s$, $F^i_j$ denotes the $i^{th}$ component of $F_j$
and `$\approx$' stands for ``up to sign".}
\vskip 3 pt
{\it Proof.} By induction over $j$. By definition of $F_0$, we have $\o_0=\text{\rm vol}_{N_{k-1}}$. Suppose
(4.41) holds for all $j\leq l<k-1$. We compute
$$
\o_{l+1}=\o_l\wedge dF^1_{l+1}\wedge dF^2_{l+1}\wedge\cdots\wedge dF^{\nu_{k-(l+1)}}_{l+1}.
$$
By Lemma 4.4, 
$$
F_{l+1}=M_{l+1}\left(\matrix \x^{N_{k-l}+1} \\ \cdot \\ \cdot \\ \x^{N_{k-(l+1)}} \\ x_{N_{k-(l+1)}+n_{k-(l+1)}+1} \\ \cdot \\ \cdot \\ x_{N_{k-(l+1)}+\nu_{k-(l+1)}}\endmatrix\right)+f_{[N_{k-l}+n_{k-l}]}\eqno(4.42)
$$
and we recall that the entries of $M_{l+1}$ are functions of rank $N_{k-l}+n_{k-l}$. By the inductive hypothesis, (4.41)
with $j=l$, we need only consider the exterior derivatives $dF^i_{l+1}$ modulo $\O_{N_{k-(l+1)}}$, where $\O_\s$ denotes the sub-bundle $\{dx,dx_1,\ldots,dx_\s\}\subseteq T^*M$.  By (4.42), we deduce that
$$
dF_{l+1}=M_{l+1}\left(\matrix d\x^{N_{k-l}+1} \\ \cdot \\ \cdot \\ d\x^{N_{k-(l+1)}} \\ dx_{N_{k-(l+1)}+n_{k-(l+1)}+1} \\ \cdot \\ \cdot \\ dx_{N_{k-(l+1)}+\nu_{k-(l+1)}}\endmatrix\right)\mod \O_{N_{k-(l+1)}}.\eqno(4.43)
$$
By Lemma 4.3, for $h$ in the range $N_{k-l}+1\leq h\leq N_{k-(l+1)}$, $\xi^h$ has rank $N_{k-(l+1)}+n_{k-(l+1)}$ and we therefore compute that
$$
d\x^h\equiv\sum_{\a=N_{k-(l+1)}+1}^{N_{k-(l+1)}+n_{k-(l+1)}}\frac{\partial \xi^h}{\partial x_\a}dx_\a\mod \O_{N_{k-(l+1)}}. 
$$
From this and (4.43), we calculate that
$$
dF_{l+1}=M_{l+1}\widehat{\G}_{l+1}\left(\matrix dx_{N_{k-(l+1)}+1} \\ dx_{N_{k-(l+1)}+2} \\ \cdot \\ \cdot
                     \\ dx_{N_{k-(l+2)}}\endmatrix\right)\mod\O_{N_{k-(l+1)}}.
$$
Consequently,
$$
\o_{l+1}=\o_l\wedge \det M_{l+1}\det\widehat{\G}_{l+1}dx_{N_{k-(l+1)}+1}\wedge\cdots
                  \wedge dx_{N_{k-(l+2)}}.\eqno(4.44)
$$
Recalling the definition of $M_l$ and that $M_1=I_{\nu_1}$, an easy induction verifies that
$$
\det M_l=\det\G_1\det\G_2\cdots\det\G_{l-1},\ 2\leq l\leq k-1.\eqno(4.45)
$$
Since $\det\widehat{\G}_{l+1}=\det\G_{l+1}$, the inductive hypothesis, together with (4.44) and (4.45) shows that 
$$
\o_{l+1}\approx (\det\G_1)^{l+1}(\det\G_2^l)\cdots (\det\G_1)^2(\det\G_{l+1})~\text{\rm vol}_{N_{k-{(l+2)}}}
$$
which proves Lemma 4.4. \qed

\vskip 5 pt

Thus we have shown that $\o_{k-1}=\Phi\hskip 0.75 pt\text{\rm vol}_{N_0}$ where, by Lemma 4.3, $\Phi$ is a generically nonzero function. Consequently,
$$
\o_{k-1}\wedge dx_{N_0+1}\wedge\cdots\wedge dx_{N_0+\nu_1}\neq 0,
$$
as we wanted to show. By their definition, the components of $F_l$ are contact coordinates. Moreover, $\CV$ is locally equivalent to the partial prolongation $\Cal C\langle\text{\rm deccel}(\CV)\rangle$. Theorem 4.1 is now proved. \qed 

\vskip 5 pt

{\smc Remark 4.1.} Theorem 4.1 reduces to the classical Goursat normal form, characterising the contact distribution on the jet bundle $J^k(\Bbb R,\Bbb R)$, in the special case when the Goursat bundle
$\CV$ has decceleration vector $\text{\rm deccel}(\CV)=\langle 0,0,0,\ldots,0,1\rangle$, where $k-1$ zeros precede the final entry, one. 
\vskip 5 pt
We conclude by presenting some simple illustrative examples. We preface this by mentioning that we have thus far left open the question of the most {\it efficient} means of constructing contact coordinates for a given partial prolongation. The proof of Theorem 4.1 is extravagant with respect to the number of integrations that are carried out in order to construct a coordinate system adapted to the natural filtration induced by $\CV$. This number can be considerably reduced, making for an efficient algorithm. However, to keep this report within bounds we will not go on to discuss this important question here. Thus, in these examples we shall be somewhat informal, being content simply, to illustrate the various constructions that we've encountered in the proof of the main result, Theorem 4.1. 
\vskip 5 pt
\noindent {\it Example 4.1.}  Consider the distribution 
$$
\CV=\Big\{\frac{1+x_2x_6}{x_6}\P {x_3}+(1+x_3)\P {x_4}+\P {x_5}+\frac{1}{x_1}\P {x_6},
\P {x_5}-x_1^2\P {x_1},\P {x_2}\Big\}
$$
defined on a generic subset of $\Bbb R^6$.
The derived type of $\CV$ is
$$
[[3,0],[5,3],[6,6]].\eqno(4.46)
$$
A calculation shows that 
$$
\ch\CV^{(1)}_0=\{\P {x_5}-x_1^2\P {x_1},\P {x_2}\},\ch\CV^{(1)}=\{\P {x_5}-x_1^2\P {x_1},\P {x_2},\P {x_6}\}\eqno(4.47) 
$$
From (4.47), and (4.46) it is easy to verify, using Proposition 3.1, that $\CV$ has the derived type of the partial prolongation $\Cal C\langle 1,1\rangle$.
Since $\D_2=1$, we have verified that $\CV$ determines a Goursat bundle of this type. By Theorem 4.1, $\CV$ is locally equivalent to the partial prolongation $\Cal C\langle 1,1\rangle$ on a generic subset $M\subset\Bbb R^6$. This settles the recognition problem for $\CV$.

Going further to compute an equivalence, the proof of Theorem 4.1, shows that this may achieved by constructing the locally defined bundle $\Pi^2$
since $\D_2=1$, in this case. To do this we must compute the invariants of $\ch\CV^{(1)}$
and choose, from any one of these, one that will be taken to be the independent variable, $x$. The invariants in question are spanned by $x_3,x_4$ and $x_5-1/x_1$. If we take $x=x_5-1/x_1$, then we find that
$$
\Pi^2=\{\P {x_5}-x_1^2\P {x_1}, \P {x_2}, \P {x_6}, \P {x_3}\}.
$$
The invariants of $\Pi^2$ are spanned by $x_5-1/x_1$ and $x_4$. From (4.47) and the invariants of $\Pi^2$, we  
deduce that $z^1=x_6$ is the coordinate that is differentiated {\it once}, while $z^2=x_4$ is the coordinate that is differentiated {\it twice}. The total differential operator is any element $X\in \CV$ that satisfies $X(x_5-1/x_1)=1$.
A convenient choice is 
$$
X=\frac{1+x_2x_6}{x_6}\P {x_3}+(1+x_3)\P {x_4}+\P {x_5}+\frac{1}{x_1}\P {x_6}
$$
from which we compute the functions
$$
z^1_1=Xz_1,\ z^2_1=Xz^2,\ z^2_2=Xz^2_1.
$$
Theorem 4.1 asserts that the functions $x,z^1,z^1_1,z^2,z^2_1,z^2_2$ are generically independent and are therefore contact coordinates for $\CV$. Indeed, we find by explicit computation that the map 
$$
\phi\:M\to J^{(1,1)}(\Bbb R,\Bbb R^2)
$$
defined by 
$$
x=x_5-1/x_1,z^1=x_6,z^1_1=1/x_1,z^2=x_4,z^2_1=1+x_3,z^2_2=\frac{1+x_2x_6}{x_6}
$$
satisfies 
$$
\phi_*\CV=\{\P x+z^1_1\P {z^1}+z^2_1\P {z^2}+z^2_2\P {z^2_1}, \P {z^1_1}, \P {z^2_2}\}
$$
as we wanted.

\vskip 5 pt

{\it Example 4.2.} Some of the contructions we've encountered in this paper are nicely illustrated by the toy control system in 3 states $x_1,x_2,x_3$ and 2 controls $u_1,u_2$, due to R. Marino (see [5]).
$$
\aligned
\dot{x}_1&=x_2+u_2x_3,
\dot{x}_2=x_3+u_2x_1,\cr
\dot{x}_3&=u_1,
\dot{x}_4=u_2.
\endaligned
$$
The corresponding vector field distribution is
$$
\CV=\big\{\P t+(x_2+u_2x_3)\P {x_1}+(x_3+u_2x_1)\P {x_2}+u_1\P {x_3}+u_2\P {x_4}, \P {u_1},\P {u_2}\big\}\eqno(4.48)
$$
whose derived type is
$$
[[3,0],[5,2],[7,7]].
$$
Thus, the derived length is 2 and  easy calculation shows that
$$
\ch\CV^{(1)}_0=\ch\CV^{(1)}=\{\P {u_1},\P {u_2}\}.
$$
From this we verify, by Proposition 3.1, that $\CV$ has the derived type of the partial prolongation $\Cal C\langle 0,2\rangle$. 
Since $\D_k=\D_2=2>1$, we must next compute the singular variety of 
$$
\aligned
\widehat{\CV}^{(1)}:=\CV^{(1)}/\ch\CV^{(1)}=\{\P t+(x_2+&u_2x_3)\P {x_1}+(x_3+u_2x_1)\P {x_2},\cr
&\ \ \ \ \P {x_3}, \P {x_4}\}=\{M_1,M_2,M_3\}.
\endaligned
$$ 
The polar matrix of a point $E=[a_1M_1+a_2M_2+a_3M_3]$ $\in\Bbb P\widehat{\CV}^{(1)}$ is 
$$
\sigma(E)=\left(\matrix -a_2 & a_1-a_3/x_1 & a_2/x_2\\
               -a_3 & -a_3x_2/x_1 & a_1+a_2x_2/x_1\endmatrix\right).
$$
We compute that the singular variety of $\s(E)$ is
$$
\widehat{\S}_1=\Bbb P\{x_1M_2-x_2M_1,M_1+x_1M_3\}\approx\Bbb R\Bbb P^1.
$$
Consequently, the resolvent bundle is
$$
\Cal R_{\widehat{\S}_1}(\CV^{(1)})=\{x_1M_2-x_2M_1,M_1+x_1M_3, \P {u_1}, \P {u_2}\},
$$
which is integrable. Thus, we've shown that $\CV$ determines a Goursat bundle of type $\langle 0,2\rangle$. By Theorem 4.1, $\CV$ is locally equivalent to the partial prolongation $\Cal C\langle 0,2\rangle$. That is, $\CV\approx\Cal C^{(2)}_2$, the contact distribution on jet bundle 
$J^2(\Bbb R,\Bbb R^2)$.

To construct an equivalence to $\Cal C^{(2)}_2$, one may follow the proof of Theorem 4.1 and compute a complete set of invariants of the resolvent bundle from which one builds a coordinate system adapted to the filtration induced by $\CV$.    
We refrain from writing out the details of this and merely remark that
the transformation so constructed is an equivalence but not a {\it static feedback} equivalence. That is, an equivalence that preserves the form of the control system and the physical distinction between time $t$, states $x_i$ and controls 
$u_\a$ (see [3,5] for details). Indeed, it is not difficult to prove that there is {\it no} static feedback equivalence between the Marino distribution (4.48) and any partial prolongation of $C^{(1)}_q$.

However, it is known that a certain {\it Cartan prolongation} (see, for example, [5]) of (4.48) {\it is} static feedback equivalent to a partial prolongation of $\Cal C^{(1)}_q$ for some $q>1$. We will consider this example as our final illustration of Theorem 4.1. 

A certain Cartan prolongation of the Marino system can be taken to be 
$$
\aligned
\dot{x}_1&=x_2+v_1x_3,\dot{x}_2=x_3+v_1x_1,\dot{x}_3=u_1,\cr
\dot{x}_4&=v_1,\dot{v}_1=v_2,\dot{v}_2=v_3,\dot{v}_3=u_2,\cr
\endaligned
$$
obtained from the Marino system by ``differentiating $u_2$ three times" (and making a slight change of notation). Here, the new state variables are $x_1,x_2,x_3,v_1,v_2,v_3$. The new control variables are $u_1$ and $u_2$. The associated distribution $\text{\rm pr}\CV$ of the Cartan prolonged Marino system is 
$$
\aligned
\text{\rm pr}\CV=\big\{\P t+(x_2+v_1x_3&)\P {x_1}+(x_3+v_1x_1)\P {x_2}+u_1\P {x_3}+\cr
&v_1\P {x_4}+v_2\P {v_1}+v_3\P {v_2}+u_2\P {v_3}, \P {u_1},\P {u_2}\big\} 
\endaligned\eqno(4.49)    
$$
whose derived type is
$$
[[3,0],[5,2],[7,4],[9,7],[10,10]]\eqno(4.50)
$$
and consequently the decceleration vector is 
$$
\text{deccel}(\text{\rm pr}\CV)=\langle 0,0,1,1\rangle. 
$$
We compute that $\text{\rm pr}\CV$ induces the filtration
$$
\aligned
\ch(\text{\rm pr}\CV)^{(1)}\subset\ch(\text{\rm pr}\CV)^{(2)}\subset&\ch(\text{\rm pr}\CV)^{(3)}_2
\subset\ch(\text{\rm pr}\CV)^{(3)}\cr
&\ \ \ \ \ \ \ \ \ \ \ \ \ \ \ \ \ \ \ \ \ \ \ \ \ \ \ \ \ \ \ \ \ \ \subset\Pi^4
\endaligned\eqno(4.51)
$$
where
$$
\aligned
&\ch(\text{\rm pr}\CV)^{(1)}=\{\P {u_1}, \P {u_2}\},\ch(\text{\rm pr}\CV)^{(2)}=\{\P {u_1},\P {u_2},\P {x_3},\P {v_3}\},\cr
&\ch(\text{\rm pr}\CV)^{(3)}_2=\ch(\text{\rm pr}\CV)^{(2)}\oplus\{\P {v_2},v_1\P {x_1}+\P {x_2}\},\cr
&\ch(\text{\rm pr}\CV)^{(3)}=\ch(\text{\rm pr}\CV)^{(3)}_2\oplus\{(1-v_2)\P {x_1}+v_1^2\P {x_2}\},\cr
&\Pi^4=\ch(\text{\rm pr}\CV)^{(3)}\oplus\{\P {v_1}\}
\endaligned
$$
and
$$
(\text{\rm pr}\CV)^{(3)}/\ch(\text{\rm pr}\CV)^{(3)}=\{\P t+v_1\P {x_4},\P {v_1}\}  
$$
The above data verifies, via Proposition 3.1, that $\text{\rm pr}\CV$ has the derived type of the partial prolongation 
$\Cal C\langle 0,0,1,1\rangle$ and that $\ch(\text{\rm pr}\CV)^{(3)}_2$ is integrable. Consequently, since $\D_k=\D_4=1$,
we have verified that $\text{\rm pr}\CV$ is a Goursat bundle on $\Bbb R^{10}$ of type $\tau=\langle 0,0,1,1\rangle$. By Theorem 4.1, the Cartan prolonged Marino distribution $\text{\rm pr}\CV$ is locally equivalent to the contact distribution 
$\Cal C\langle 0,0,1,1\rangle$. But is it static feedback equivalent?

Taking coordinates adapted to filtration (4.51) reveals that $t$ may be taken to be the independent variable $x$, while
$x_4$ may be taken to be the coordinate that must be differentiated 4 times. Finally $x_1-v_1x_2$ may be taken to be the variable that is differentiated 3 times. Consequently, the total derivative operator is
$$
\aligned
X=\P t+(x_2+v_1x_3)\P {x_1}&+(x_3+v_1x_1)\P {x_2}+u_1\P {x_3}+\cr
&\ \ \ \ \ \ \ v_1\P {x_4}+v_2\P {v_1}+v_3\P {v_2}+u_2\P {v_3}
\endaligned
$$
and we construct contact coordinates as in Theorem 4.1:
$$
\aligned
&x=t,z^1=x_1-v_1x_2,z^1_1=Xz^1,z^1_2=Xz^1_1,z^1_3=Xz^1_2,\cr
&z^2=x_4,z^2_1=Xz^2,z^2_2=Xz^2_1,z^2_3=Xz^2_2,z^2_4=Xz^2_3
\endaligned\eqno(4.52)
$$
leading to the equivalence
$$
\psi\:\Bbb R^{10}\to J^{(3,4)}(\Bbb R,\Bbb R^2),
$$
defined by
$$
\aligned
&x=t,\cr
&z^1=x_1-v_1x_2,\cr
&z^1_1=x_2-x_1v_1^2-v_2x_2,\cr
&z^1_2=-v_1^2x_2-v_1^3x_3+x_3-x_3v_2+x_1v_1-3v_2x_1v_1-v_3x_2,\cr
&z^1_3= v_1x_2-5x_2v_2v_1-6x_3v_1^2v_2-2x_3v_3-\cr
& \  \ \ \ \ \ \ \ \ x_1v_1^3-4x_1v_1v_3-u_1v_1^3+u_1-u_1v_2+v_2x_1-3v_2^2x_1-u_2x_2, \cr
&z^2=x_4 ,
z^2_1=v_1,
z^2_2=v_2,
z^2_3=v_3,
z^2_4=u_2.
\endaligned
$$
By Theorem 4.1, $\psi$ is a local diffeomorphism. Indeed the determinant of the Jacobian matrix of $\psi$
is  $-(-1+v_2+v_1^3)^3$. By its construction, $\psi$ is an equivalence between $\text{\rm pr}\CV$, the prolonged
Marino distribution (4.49) and the partial prolongation $\Cal C\langle 0,0,1,1\rangle$, of $\Cal C^{(1)}_2$.
From the form of $\psi$, we observe that it's a static feedback equivalence.  This then reveals the fact, well known for this example, that while the Marino distribution (4.48) is not static feedback equivalent to  a partial prolongation of $\Cal C^{(1)}_q$ for some $q\geq 1$, nevertheless, some Cartan prolongation of it does admit such an equivalence.

\vskip 20 pt

\centerline{\it References}
%\item{[1]} D. Bao, S-S. Chern, Z. Shen (Eds), Finsler Geometry, Joint summer research conference on Finsler geometry,
%{\it Contemporary Mathematics Vol. 196}, American Mathematical Society, 1996
%\item{[2]} P. Brunovsk\'y, A classification of linear controllable systems, {\it Kybernetika}, \v Cislo 3, Ro\v cnik 6, %(1970), 173-188
\item{[1]} R.L. Bryant, {\it Some Aspects of the Local and Global Theory of Pfaffian Systems}, Ph.D thesis, University of North Carolina, Chapel Hill, 1979
\item{[2]} R.L. Bryant, S-S. Chern, R.B. Gardner, H. Goldschmidt, P. Griffiths, {\it  Exterior Differential Systems},
MSRI Publications, \allowbreak Springer-Verlag, New York, 1991
%\item{\bf 5.} \'E. Cartan, Sur l'\'equivalence absolue de certaines syst\`emes  \newline d'\'equations diff\'erentielles et certaines syst\`emes de courbes, {\it Bulletin de la Soci\'et\'e Math\'ematiques de France}, 42, (1914), 12-48
%\item{\bf 6.} \'E. Cartan, Sur l'int\'egration de certains syst\'emes ind\'etermin\'es d'\'equations diff\'erentielles, 
{%\it J. f\"ur Reine und Angew. Math.}, 145, (1915), 86-91  
%\item{\bf 7.} R.B. Gardner, Invariants of Pfaffian systems, {\it Trans. Amer. Math. Soc.}, 126(3), (167), 514-533
\item{[3]} R.B. Gardner, W.F. Shadwick, The GS algorithm for exact linearization to Brunovsk\'y normal form, {\it IEEE Trans. Automat. Control}, 37(2) (1992), 224-230
%\item{\bf 9.} A. Giaro, A. Kumpera, C. Ruiz, Sur la lecture correcte d'un r\'esult d'\'Elie Cartan, 
%{\it C.R. Acad. Sc. Fr.}, 287, S'erie A, (1978), 241-244

%\item{\bf 10.} E. Goursat, {\it Le\c cons sur le probl\`eme de Pfaff}, Hermann, Paris, 1923
%\item{\bf 11.} P. Libermann, Sur le probl\`eme d'\'equivalence des syst\`emes de Pfaff non compl\`ement integrables, in {\it Differential Geometry}, Publ. \newline Math. Univ. Paris VII, Volume 3, 73-109, Universit\'e de Paris VII, 1978  
%\item{\bf 12.} J. Martinet, Classes charact\'eristiques des syst\`emes de Pfaff,  \newline{\it Journ\'ees Exotiques}, Lille, 1973
\item{[4]} W. Respondek, W. Pasillas-Lepine, %Extended Goursat normal \newline form:a geometric characterisation, in Nonlinear Control in the Year 2000, {\it Lecture Notes in Control and Information Sciences}, Volume 259, 323-338, Springer-Verlag, 2001; 
Contact systems and corank 1 involutive subdisributions, {\it Acta Appl. Math.}, {\bf 69} (2001), 105-128
\item{[5]} W. Sluis, Absolute equivalence and its applications to control theory, 33rd IEEE Conference on Decision and Control, Lake Buena Vista, Florida, December, 1994 available at\newline 
{\tt www.cds.caltech.edu/\ $\tilde{}$ murray/archive/eds-dec96/cdc94-workshop.html}  
\item{[6]} O. Stormark, {\it Lie's Structural Approach to PDE Systems}, Encyclopedia of Mathematics and its Applications, vol {\bf 80}, Cambridge University Press, 2000
%\item{\bf 14.} D. Tilbury, R.M. Murray, S. Shankar, Trajectory generation for the $N$-trailer problem using the Goursat normal form, {\it IEEE Trans. Automat. Control}, 40(5), (1995), 802-819
%\item{\bf 15.} P.J. Vassiliou, Local equivalence for vector field systems, {\it Bull. Austral. Math. Soc.}, 42(2), (1990), 215-229
\item{[7]} E. von Weber, Zur invariantentheorie der systeme Pfaff'scher gleichungen, \newline {\it Berichte Verhandlungen
der Koniglich Sachsischen Gesellshaft der Wissen-\newline shaften Mathematisch Physikalische Klasse, Leipzig}, 50: 207-229, (1898)
\item{[8]} K. Yamaguchi, Contact geometry of higher order, {\it Japan J.} 
 {\it Math.}, 8, (1982), 109-176; Geometrization of jet bundles, {\it Hokkaido Math. J.}, 12, (1983), 27-40

\vskip 10pt
{\smc Acknowledgments}: It is a pleasure to thank George Wilkens for his reading of a previous version of this paper and for his valuable input. The computations in this paper were carried out using the general purpose differential geometry package
{\smc Vessiot} developed at Utah State University by Ian M. Anderson and colleagues. It is a pleasure to thank the School of Mathematics and Statistics at Utah State and Ian in particular for hospitality in 2002 when a portion of the research reported here was carried out.    
\vskip 10 pt
%\address School of Mathematics and Statistics,
%University of Canberra, Australian Capital Territory, AUSTRALIA., 2601\endaddress

\enddocument
\bye